\pdfoutput=1
\RequirePackage{ifpdf}
\ifpdf 
\documentclass[pdftex]{sigma}
\else
\documentclass{sigma}
\fi

\numberwithin{equation}{section}
\newtheorem{Theorem}{Theorem}[section]
\newtheorem{Corollary}[Theorem]{Corollary}
\newtheorem{Lemma}[Theorem]{Lemma}
\newtheorem{Proposition}[Theorem]{Proposition}
{\theoremstyle{definition}
\newtheorem{Definition}[Theorem]{Definition}
\newtheorem{Remark}[Theorem]{Remark}
}

\DeclareMathOperator{\bdd}{bdd}
\DeclareMathOperator{\ch}{ch}
\DeclareMathOperator{\conv}{conv}
\DeclareMathOperator{\End}{End}
\DeclareMathOperator{\Ext}{Ext}
\DeclareMathOperator{\ev}{ev}
\DeclareMathOperator{\gr}{gr}
\DeclareMathOperator{\Hom}{Hom}
\DeclareMathOperator{\id}{id}
\DeclareMathOperator{\Ob}{Ob}
\DeclareMathOperator{\soc}{soc}
\DeclareMathOperator{\wt}{wt}

\begin{document}

\allowdisplaybreaks

\renewcommand{\thefootnote}{$\star$}

\renewcommand{\PaperNumber}{030}

\FirstPageHeading

\ShortArticleName{Tilting Modules in Truncated Categories}

\ArticleName{Tilting Modules in Truncated Categories\footnote{This paper is a~contribution to the Special Issue on New
Directions in Lie Theory.
The full collection is available at \href{http://www.emis.de/journals/SIGMA/LieTheory2014.html}
{http://www.emis.de/journals/SIGMA/LieTheory2014.html}}}

\Author{Matthew BENNETT~$^\dag$ and Angelo BIANCHI~$^\ddag$}

\AuthorNameForHeading{M.~Bennett and A.~Bianchi}

\Address{$^\dag$~Department of Mathematics, State University of Campinas, Brazil}
\EmailD{\href{mailto:mbenn002@gmail.com}{mbenn002@gmail.com}}

\Address{$^\ddag$~Institute of Science and Technology, Federal University of S\~ao Paulo, Brazil}
\EmailD{\href{mailto:acbianchi@unifesp.br}{acbianchi@unifesp.br}}

\ArticleDates{Received September 05, 2013, in f\/inal form March 17, 2014; Published online March 26, 2014}

\Abstract{We begin the study of a~tilting theory in certain truncated categories of mo\-du\-les $\mathcal G(\Gamma)$ for the
current Lie algebra associated to a~f\/inite-dimensional complex simple Lie algebra, where $\Gamma = P^+ \times J$, $J$ is
an interval in $\mathbb Z$, and $P^+$ is the set of dominant integral weights of the simple Lie algebra.
We use this to put a~tilting theory on the category $\mathcal G(\Gamma')$ where $\Gamma' = P' \times J$, where
$P'\subseteq P^+$ is saturated.
Under certain natural conditions on $\Gamma'$, we note that $\mathcal G(\Gamma')$ admits full tilting modules.}


\Classification{17B70; 17B65; 17B10; 17B55}

\renewcommand{\thefootnote}{\arabic{footnote}}
\setcounter{footnote}{0}

\section{Introduction}

Associated to any f\/inite-dimensional complex simple Lie algebra $\mathfrak g$ is its current algebra $\mathfrak g[t]$.
The current algebra is just the Lie algebra of polynomial maps from $\mathbb C \rightarrow \mathfrak g$ and can be
identif\/ied with the space $\mathfrak g \otimes \mathbb C[t]$ with the obvious commutator.
The study of the representation theory of current algebras was largely motivated by its relationship to the
representation theory of af\/f\/ine and quantum af\/f\/ine algebras associated to $\mathfrak g$.
However, it is also now of independent interest since the current algebra has connections with problems arising in
mathematical physics, for instance the $X = M$ conjectures, see~\cite{AK,FK,naoi1}.
Also, the current algebra, and many of its modules, admits a~natural grading by the integers, and this grading gives
rise to interesting combinatorics.
For example,~\cite{KNloewy} relates certain graded characters to the Poincar\'e polynomials of quiver varieties.

Let $P^+$ be the set of dominant integral weights of $\mathfrak g$, $\Lambda=P^+ \times \mathbb Z$, and
$\widehat{\mathcal G}$ the category of $\mathbb Z$-graded modules for $\mathfrak g[t]$ with the restriction that the
graded pieces are f\/inite-dimensional.
Also, let $\mathcal G$ be the full subcategory of $\widehat{\mathcal G}$ consisting of modules whose grades are bounded
above.
Then $\Lambda$ indexes the simple modules in $\widehat{\mathcal G}$.
In this paper we are interested in studying Serre subcategories $\widehat{\mathcal G}(\Gamma)$ where $\Gamma \subset
\Lambda$ is of the form $P'\times J$ where $J\subset \mathbb Z$ is a~(possibly inf\/inite) interval and $P'\subset P^+$ is
closed with respect to a~natural partial order.
In particular, we study the tilting theories in these categories.
This generalized the work of~\cite{bc}, where $\Gamma$ was taken to be all of $\Lambda$.

The category $\widehat{\mathcal G}(\Gamma)$ contains the projective cover and injective envelope of its simple objects.
Given a~partial order on the set $\Gamma$, we can def\/ine the standard and costandard objects, as in~\cite{Donkin}.
The majority of the paper is concerned with a~particular order, in which case the standard objects $\Delta(\lambda,
r)(\Gamma)$ are quotients of the f\/inite-dimensional local Weyl modules, and the costandard objects $\nabla(\lambda,
r)(\Gamma)$ are submodules of (appropriately def\/ined) duals of the inf\/inite-dimensional global Weyl modules.
We recall (see, for example,~\cite{Mathieu}) that a~module~$T$ is called tilting if~$T$ admits a~f\/iltration by standard
modules and a~f\/iltration by costandard modules.
In our case, both sets of objects have been extensively studied (see~\cite{CL,FKtw,FoL,naoi2} for the local Weyl
modules, and~\cite{CFK,CPweyl}, for the global Weyl modules).
Both families of modules live in a~subcategory $\mathcal G_{\bdd}(\Gamma)$ consisting of objects whose
weights are in a~f\/inite union of cones (as in $\mathcal O$) and whose grades are bounded above.
The main goal of this paper is to construct another family of modules indexed by $\Gamma$ and which are in $\mathcal
G_{\bdd}(\Gamma)$.
These modules are denoted by $T(\lambda, r)(\Gamma)$, and admit an inf\/inite f\/iltration whose quotients are of the form
$\Delta(\mu,s)(\Gamma)$, for $(\mu,s)\in \Gamma$.
They also satisfy the homological property that $\Ext^1_{\mathcal G}(\Delta(\mu,s)(\Gamma), T(\lambda,
r)(\Gamma))=0$ for all $(\mu, s)\in \Gamma$.
We use the following theorem to prove that this homological property is equivalent to having
a~$\nabla(\Gamma)$-f\/iltration, proving that the $T(\lambda, r)(\Gamma)$ are tilting.
The theorem was proved in~\cite{bcm,bbcklk}, and~\cite{ci} for $\mathfrak sl_2[t]$, $\mathfrak sl_{n+1}[t]$ and
general $\mathfrak g[t]$ respectively.

\begin{Theorem}
Let $P(\lambda, r)$ denote the projective cover of the simple module $V(\lambda, r)$.
Then $P(\lambda,r)$ admits a~filtration by global Weyl modules, and we have an equality of filtration multiplicities
$[P(\lambda, r): W(\mu,s)] = [\Delta(\mu,r): V(\lambda,s)]$, where $\Delta(\mu,r)$ is the local Weyl module.
\end{Theorem}

The following is the main result of this paper.

\begin{Theorem}\quad
\begin{enumerate}\itemsep=0pt
\item[$1.$]
Given $(\lambda, r)\in\Gamma$, there exists an indecomposable module $T(\lambda,r)(\Gamma)\in\Ob\mathcal
G_{\bdd}(\Gamma)$ which admits a~$\Delta(\Gamma)$-filtration and a~$\nabla(\Gamma)$-filtration.
Further,
\begin{gather*}
T(\lambda,r)(\Gamma)[s]_\lambda=0
\qquad
\text{if}
\quad
s>r,
\quad
T(\lambda,r)(\Gamma)[r]_\lambda=1,
\quad
\wt T(\lambda,r)(\Gamma)\subset\conv W\lambda,
\end{gather*}
and $T(\lambda,r)(\Gamma)\cong T(\mu,s)(\Gamma)$ if and only if $(\lambda,r)=(\mu,s)$.
\item[$2.$]
Moreover, any indecomposable tilting module in $\mathcal G_{\bdd}(\Gamma)$ is isomorphic to
$T(\lambda,r)(\Gamma)$ for some $(\lambda,r)\in \Gamma$, and any tilting module in $\mathcal
G_{\bdd}(\Gamma)$ is isomorphic to a~direct sum of indecomposable tilting modules.
\end{enumerate}
\end{Theorem}

The majority of the paper is devoted to the case where $\Gamma=P^+ \times J$.
It is easy to see from the construction that the module $T(\lambda, r)(\Gamma)$ has its weights bounded above by
$\lambda$.
It follows that if we let $P'\subset P^+$ be saturated (downwardly closed with respect to the normal partial order on
weights), and set $\Gamma'=P'\times J$, then $T(\lambda,r)(\Gamma')=T(\lambda, r)(\Gamma)$.

We use the convention that $\Delta(\lambda, r)(\Lambda)$ is simply written $\Delta(\lambda, r)$, and similarly for other
objects.
Keeping $\Gamma=P^+ \times J$, there is a~natural functor taking $M\in \Ob \mathcal G$ to
$M^\Gamma\in\Ob\mathcal G(\Gamma)$.
For $(\lambda, r)\in \Gamma$ this functor preserves many objects, and in particular we have $\Delta(\lambda,
r)^\Gamma=\Delta(\lambda, r)(\Gamma)$ and $\nabla(\lambda, r)^\Gamma=\nabla(\lambda, r)(\Gamma)$.
So it is natural to ask if $T(\lambda, r)^\Gamma=T(\lambda, r)(\Gamma)$.
The answer is ``no", and is a~result of the following phenomena: for $(\mu,s)\notin \Gamma$, the module
$\nabla(\mu,s)^\Gamma$ is not in general zero, and does not correspond to any simple module.
Hence $\nabla(\mu,s)^\Gamma$ can not be considered costandard.
So the modules $T(\lambda, r)(\Gamma)$ must be studied independently.

Another purpose of this paper is the following.
In~\cite{bc}, the tilting modules $T(\lambda, r)$ are constructed for all $(\lambda, r)\in \Lambda$.
It is normal to then consider the module $T=\bigoplus_{(\lambda,r)\in \Lambda} T(\lambda, r)$, the algebra
$A=\End T$, and use several functors to f\/ind equivalences of categories.
However, it is not hard to see that if $T$ is def\/ined in this way, then $T$ fails to have f\/inite-dimensional graded
components, and hence $T\notin\Ob\mathcal G$.
One of the purposes of this paper is to f\/ind Serre subcategories with index sets $\Gamma$ such that
$T(\Gamma)=\bigoplus_{(\lambda,r)\in\Gamma}T(\lambda,r)(\Gamma)\in\Ob\mathcal G(\Gamma)$.
It is not hard to see that (except for the degenerate case where $\Gamma=\{0\}\times J$) a~necessary and suf\/f\/icient pair
of conditions on $\Gamma$ is that $P'$ be f\/inite and $J$ have an upper bound.
It is natural to study the algebra $\End T(\Gamma)$ in the case that $T(\Gamma ) \in \mathcal G(\Gamma)$,
and this will be pursued elsewhere.
We also note that in the case that $\Gamma$ is f\/inite then $\End T(\Gamma)$ is a~f\/inite-dimensional
associative algebra.

We end the paper by considering other partial orders which can be used on $\Gamma\subset\Lambda$.
In particular, we consider partial orders induced by the so-called covering relations.
One tends to get trivial tilting theories in these cases (one of the standard-costandard modules is simple, and the
other is projective or injective), but the partial orders are natural for other reasons, and we include their study for
completeness.
One of the reasons to study these other subcategories is that one can obtain directed categories as in~\cite{cg:hwcat}
(in the sense of~\cite{CPS}).

The paper is organized as follows.
In Section~\ref{section1}
we establish notation and recall some basic results on the f\/inite-dimensional representations of
a~f\/inite-dimensional simple Lie algebra.
In Section~\ref{section2}
we introduce several important categories of modules for the current algebra.
We also introduce some important objects, including the local and global Weyl modules.
In Section~\ref{section3} we state the main results of the paper and establish some homological results.
Section~\ref{section4} is devoted to constructing the modules $T(\lambda, r)(\Gamma)$ and establishing their properties.
Finally, in Section~\ref{section5}, we consider the tilting theories which arise when considering partial orders on $\Lambda$ which
are induced by covering relations.

We also provide for the reader's convenience a~brief index of the notation which is used repeatedly in this paper.

\section{Preliminaries}\label{section1}

\subsection{Simple Lie algebras and current algebras}\label{section1.1}

We f\/ix $\mathfrak g$,  a~complex simple f\/inite-dimensional Lie algebra, and let $\mathfrak h \subset \mathfrak g$ be a~f\/ixed Cartan subalgebra.  Denote by $\{\alpha_i: i\in I\}$ a set of simple roots of $\mathfrak g $ with respect to the Cartan subalgebra $\mathfrak h$, where $I=\{1,\dots,\dim\mathfrak h\}$.   Let $R\subset\mathfrak h^*$ be the corresponding set of roots, $R^+$ the positive roots, $P^+$ the dominant integral weights, and $Q^+$ the positive root lattice.  By $\theta$ we denote the highest root.  Given $\lambda,\mu\in\mathfrak h^*$, we say that  $\lambda \ge \mu$ if and only if $ \lambda-\mu\in Q^+$.  The Weyl group of $\mathfrak g$ is the subgroup $W\subset \operatorname{Aut}(\mathfrak h^*)$ generated by the simple ref\/lections~$s_i$, and we let~$w_\circ$ denote the unique longest element of~$W$.   For $\alpha\in R$ we write $\mathfrak g_\alpha$ for the corresponding root space.  Then the subspaces $\mathfrak n^\pm=\bigoplus_{\alpha\in R^+}\mathfrak g_{\pm\alpha},$ form Lie subalgebras of $\mathfrak g$.  We f\/ix a Chevalley basis $\{ x^\pm_\alpha, \, h_i \,|\, \alpha\in R^+, \, i\in I \}$ of $\mathfrak g$, and for each $\alpha\in R^+$ we set $h_\alpha=[x_\alpha,x_{-\alpha}]$. Note that $h_{\alpha_i}=h_i$, $i\in I$, and we let  $\omega_i=h_i^*\in P^+$.

For any Lie algebra $\mathfrak a$ we can construct another Lie algebra $\mathfrak a[t]=\mathfrak a \otimes \mathbb C[t]$, with bracket given by $[x\otimes t^r, y\otimes t^s]=[x,y]\otimes t^{r+s}$, which is the current algebra associated to $\mathfrak a$.  Set $\mathfrak a[t]_+=\mathfrak a \otimes t \mathbb C[t]$.  Then $\mathfrak a[t]$ and $\mathfrak a[t]_+$ are $\mathbb Z_+$-graded Lie algebras, graded by powers of $t$.  If we denote by $\mathbb U(\mathfrak a)$ the universal enveloping algebra of a Lie algebra $\mathfrak a$, then $\mathbb U(\mathfrak a[t])$ and $\mathbb U(\mathfrak a[t]_+)$ inherit a natural grading by powers of $t$.  We denote by $\mathbb U(\mathfrak a[t])[k]$ the $k^{th}$-graded component.  Each graded component is a module for $\mathfrak a$ under left or right multiplication, and the adjoint action.  Supposing that  $\dim \mathfrak a < \infty$, then the graded component $\mathbb U(\mathfrak a[t])[k]$ is a free $\mathfrak a $ module (under multiplication) of f\/inite rank.

It is well-known that the universal enveloping algebra $\mathbb U(\mathfrak a)$ is a Hopf algebra. In particular, it is equipped with a comultiplication def\/ined by sending  $x\to x\otimes 1+1\otimes x$ for $x\in\mathfrak a$, and extending this assignment to be a homomorphism.  In the case where $\mathfrak a =\mathfrak b[t]$ or $\mathfrak b[t]_+$, the comultiplication is a homomorphism of graded associative algebras.  We note that if $[\mathfrak a,\mathfrak a]=\mathfrak a$ (which holds for our Lie algebra $\mathfrak g$), then as a graded associative algebra,  $\mathfrak a$ and $\mathfrak a\otimes t$ generate ${\textbf U}(\mathfrak a[t])$.

\subsection{Finite-dimensional modules}\label{section1.2}

The f\/irst category we consider is $\mathcal F(\mathfrak g)$ the category of f\/inite-dimensional modules for~$\mathfrak g$ with morphisms $\mathfrak g$-module homomorphisms.  It is well known that this is a semi-simple category, and that the simple objects are parametrized by $\lambda \in P^+$. Letting $V(\lambda)$ denote the simple module associated to $\lambda$, it is generated by a vector $v_\lambda\in V(\lambda)$ satisfying the def\/ining relations
\begin{gather*}
  \mathfrak n^+ v_\lambda=0,\qquad hv_\lambda=\lambda(h)v_\lambda,\qquad (x^-_{\alpha_i})^{\lambda(h_i)+1}v_\lambda =0,
  \end{gather*}
for all $h\in\mathfrak h$, $i\in I$.  This category admits a duality, which on simple modules is given by $V(\lambda)^*\cong_{\mathfrak g} V(-w_\circ\lambda)$.  An object $V\in\mathcal F(\mathfrak g)$ has a weight space decomposition $V=\bigoplus_{\lambda\in\mathfrak h^*} V_\lambda$ where $V_\lambda=\{v\in V : hv=\lambda(h)v,\, \forall\, h\in\mathfrak h\}$.  For any such $V$, we def\/ine the subset  $\operatorname{wt}(V)=\{\lambda\in\mathfrak h^*:V_\lambda\ne 0\}$  and def\/ine the character of $V$ to be the sum $\ch V=\sum \dim V_\lambda r^\lambda$. The following results are standard:

\begin{Lemma} Let $V\in\mathcal F(\mathfrak g)$ and $\lambda\in P^+$. Then,
\begin{enumerate}\itemsep=0pt
\item[$1)$]  $w\operatorname{wt}(V)\subset\operatorname{wt}(V)$ and $\dim V_\lambda=V_{w\lambda}$ for all $w\in W$;
\item[$2)$]  $\dim\operatorname{Hom}_{\mathfrak g}(V(\lambda), V)=\dim\{v\in V_\lambda: \mathfrak n^+ v=0\}$;
\item[$3)$]  $\operatorname{wt}(V(\lambda))\subset \lambda-Q^+$.
\end{enumerate}
\end{Lemma}

\section{The main category and its subcategories} \label{section2}

In this section we introduce the main categories of study, and present several properties and functors between them.  We will also introduce several families of modules which will play important roles.  Most of these categories and objects have been studied elsewhere (see~\cite{bbcklk,bc,cg:hwcat}).

\subsection{The main category}\label{section2.1}

We denote by $\widehat{\mathcal G}$ the category of $\mathbb Z$-graded $\mathfrak g[t]$ modules such that the graded components are f\/inite-dimensional and where morphisms are degree zero maps of $\mathfrak g[t]$-modules.  Writing $V\in\operatorname{Ob}\widehat{\mathcal G}$ as
\begin{gather*}
 V=\bigoplus_{r\in\mathbb Z} V[r],
\end{gather*}
we see that $V[r]$ is a f\/inite-dimensional $\mathfrak g$ module, while $z\otimes t^k. V[r]\subset V[r+k]$ for all $z\in\mathfrak g$, $k\in \mathbb Z_{\ge 0}$, and $r\in \mathbb Z$.  For $M\in \widehat{\mathcal G}$,  its graded character is the sum (formal, and possibly inf\/inite)
$\ch_{\operatorname{gr}} M = \sum\limits_{r\in \mathbb Z} \ch M[r]u^r$.

For $V\in \mathcal F(\mathfrak g)$ we make $V$ an object in $\widehat{\mathcal G}$, which we shall call $\ev V$, in the following way.  Set $\ev V[0]=V$ and $\ev V[r]=0$ for all $r\ne 0$.  Then necessarily we have $z\otimes t^k. v = \delta_{k,0} z.v$ for $z\in\mathfrak g$, $k\in\mathbb Z_+$, $v\in \ev V$.  It is not hard to see that this def\/ines a covariant functor $\ev:\mathcal F(\mathfrak g)\to \widehat{\mathcal G}$.  Further, for $s\in\mathbb Z$ let $\tau_s:\widehat{\mathcal G}\to\widehat{\mathcal G}$ be the grade shift functor given by
\begin{gather*}
(\tau_s V)[k]=V[k-s],\qquad \text{for all} \quad  k\in\mathbb Z, \quad V\in\operatorname{Ob}{\widehat{\mathcal G}}.
\end{gather*}
For $(\lambda,r)\in P^+\times \mathbb Z$  set $V(\lambda,r):=\tau_r(\ev (V(\lambda)))$ and $v_{\lambda,r}:=\tau_r(v_\lambda)$.

\begin{Proposition}
The isomorphism classes of simple objects in $\widehat{\mathcal G}$ are parametrized by pairs $(\lambda,r)$ and we have
\begin{gather*}
\operatorname{Hom}_{\widehat{\mathcal G}}(V(\lambda,r), V(\mu,s))=
\begin{cases} 0, & \text{if} \ \ (\lambda,r)\ne (\mu,s),\\ \mathbb C, & \text{if} \ \ (\lambda,r) = (\mu,s).
 \end{cases}
\end{gather*}
Moreover, if $V\in\operatorname{Ob}\widehat{\mathcal G}$ satisfies $V=V[n]$ for some $n\in\mathbb Z$, then $V$ is semi-simple.
\end{Proposition}

The category $\widehat{\mathcal G}$ admits a duality, where given $M$ we def\/ine $M^*\in \operatorname{Ob}\widehat{\mathcal G}$ to be the module given by
\begin{gather*}
 M^*=\bigoplus_{r\in\mathbb Z}M^*[-r] \qquad \text{and} \qquad M^*[-r]=M[r]^*
 \end{gather*}
and equipped with the usual action where
\begin{gather*}
\big(x\otimes t^r\big)m^*(v) = -m^*\big(x\otimes t^r.v\big).
\end{gather*}
We note that $M^{**}\cong M$ and that $\ch_{\operatorname{gr}}M^*=\sum\limits_{r\in\mathbb Z} \ch (M[r]^*)u^{-r}$.

Denote by $\Lambda=P^+\times\mathbb Z$ and equip $\Lambda$ with the lexicographic partial order $\le$, i.e.\
\begin{gather*}
(\mu,r) \le (\lambda,s) \qquad \Leftrightarrow \qquad \text{either} \quad \mu < \lambda \quad \text{or} \quad \mu=\lambda \quad \text{and} \quad r\le s.
\end{gather*}

\subsection{Some bounded subcategories of the main category}\label{section2.2}

We let $\mathcal G_{\le s}$ be the full subcategory of $\widehat{\mathcal G}$ whose objects $V$ satisfy $V[r]=0$ for all $r>s$.  Clearly~$\mathcal G_{\le s}$ is a full subcategory of $\mathcal G_{\le r}$ for all $s< r\in\mathbb Z$. Def\/ine $\mathcal G$ to be the full subcategory of $\widehat{\mathcal G}$ whose objects consist of those objects $V$ satisfying $V\in\operatorname{Ob}\mathcal G_{\le s}$ for some $s\in\mathbb Z$. Finally, let~${\mathcal G}_{\operatorname{bdd}}$ be the full subcategory of~$\mathcal G$ consisting of objects~$M$ satisfying the following condition: $|\wt (M)\cap P^+|<\infty$.

 Given $s\in\mathbb Z$ and $V\in \widehat{\mathcal G}$, def\/ine a submodule $V_{>s}=\bigoplus_{r>s} V[r]$ and a corresponding quotient $V_{\le s}= V/V_{>s}$.  Then it is clear that  $V_{\le s}\in\operatorname{Ob}\mathcal G_{\le s}$, and indeed this is the maximal quotient of~$V$ in~$\mathcal G_{\le s}$.  Any $f\in\operatorname{Hom}_{\widehat{\mathcal G}}( V, W)$ naturally induces a morphism $f_{\le s}\in\operatorname{Hom}_{\widehat{\mathcal G}_{\le s}}( V_{\le s}, W_{\le s})$.  The following is proved in~\cite{cg:hwcat}.

\begin{Lemma}\label{functor} 
The assignments $V\mapsto V_{\le r}$ for all $V\in\operatorname{Ob}\widehat{\mathcal G}$ and $f\mapsto f_{\le r}$ for all $f\in\operatorname{Hom}_{\widehat{\mathcal G}}(V,W)$, $V,W\in\operatorname{Ob}\widehat{\mathcal G}$, define a full, exact and essentially surjective functor from $\widehat{\mathcal G}$ to $\mathcal G_{\le r}$.
\end{Lemma}

Given $V\in\operatorname{Ob}\mathcal G$ def\/ine $[V: V(\lambda, r)]:=[V[r]:V(\lambda)]$ the multiplicity of $V(\lambda)$ in a composition series for $V[r]$ as a~$\mathfrak g$-module. For any $V\in\operatorname{Ob}\widehat{\mathcal G}$, def\/ine
\begin{gather*}
\Lambda(V)=\{(\lambda,r)\in\Lambda: [V:V(\lambda,r)]\ne 0\}.
\end{gather*}

\subsection{Projective and injective objects in the main category}\label{section2.3}

Given $(\lambda,r)\in \Lambda$, set
\begin{gather*}
P(\lambda,r)=\textbf U(\mathfrak g[t])\otimes_{\textbf U(\mathfrak g)} V(\lambda,r)
\qquad
\text{and}
\qquad
I(\lambda,r)\cong P(-w_\circ\lambda,-r)^{*}.
\end{gather*}
Clearly these are an inf\/inite-dimensional $\mathbb Z$-graded $\mathfrak g[t]$-module.
Using the PBW theorem we have an isomorphism of graded vector spaces $\textbf U(\mathfrak g[t])\cong \textbf U(\mathfrak
g[t]_+)\otimes \textbf U(\mathfrak g)$, and hence we get $P(\lambda,r)[k]=\textbf U(\mathfrak g[t]_+)[k-r]\otimes
V(\lambda,r)$, where we understand that $\textbf U(\mathfrak g[t]_+)[k-r]=0$ if $k<r$.
This shows that $P(\lambda,r)\in\Ob\widehat{\mathcal G}$ and also that $P(\lambda,r)[r] =1\otimes
V(\lambda,r)$.
Set $p_{\lambda,r}=1\otimes v_{\lambda,r}$.

\begin{Proposition}
\label{projectives}
Let $(\lambda,r)\in\Lambda$, and $s\ge r$.
\begin{enumerate}\itemsep=0pt
\item[$1.$]
$P(\lambda,r)$ is generated as a~$\mathfrak g[t]$-module by $p_{\lambda,r}$ with defining relations
\begin{gather*}
(\mathfrak n^+)p_{\lambda,r}=0,
\qquad
hp_{\lambda,r}=\lambda(h)p_{\lambda,r},
\qquad
(x^-_{\alpha_i})^{\lambda(h_i)+1}p_{\lambda,r}=0,
\end{gather*}
for all $h\in\mathfrak h$, $i\in I$.
Hence, $P(\lambda,r)$ is the projective cover in the category $\widehat{\mathcal G}$ of its unique simple quotient
$V(\lambda,r)$.
\item[$2.$]
The modules $P(\lambda,r)_{\le s}$ are projective in $\mathcal G_{\le s}$.
\item[$3.$]
Let $V\in\Ob\widehat{\mathcal G}$.
Then $\dim\Hom_{\widehat{\mathcal G}}(P(\lambda,r),V)=[V:V(\lambda,r)]$.
\item[$4.$]
Any injective object of $\mathcal G$ is also injective in $\widehat{\mathcal G}$.
\item[$5.$]
Let $(\lambda,r)\in\Lambda$.
The object $I(\lambda,r)$ is the injective envelope of $V(\lambda,r)$ in $\mathcal G$ or $\widehat{\mathcal G}$.
\end{enumerate}
\end{Proposition}

\subsection{Local and global Weyl modules}\label{section2.4}

The next two families of modules in ${\mathcal G}_{\bdd}$ we
need are the local and global Weyl modules which were originally def\/ined in~\cite{CPweyl}.

For the purposes of this paper, we shall denote the local Weyl modules by $\Delta(\lambda,r)$, $(\lambda,r)\in
P^+\times\mathbb Z$.
Thus, $\Delta(\lambda, r)$ is generated as a~$\mathfrak g[t]$-module by an element $w_{\lambda,r}$ with relations:
\begin{gather*}
\mathfrak n^+[t]w_{\lambda,r}=0,
\qquad
(x_{i}^-)^{\lambda(h_{i})+1}w_{\lambda,r}=0,
\qquad
(h\otimes t^s)w_{\lambda,r}= \delta_{s,0}\lambda(h)w_{\lambda,r},
\end{gather*}
here $i\in I$, $h\in\mathfrak h$ and $s\in\mathbb Z_+$.

Next, let $W(\lambda,r)$ be the global Weyl modules, which is $\mathfrak g[t]$-module generated as a~$\mathfrak
g[t]$-module by an element $w_{\lambda,r}$ with relations:
\begin{gather*}
\mathfrak n^+[t]w_{\lambda,r}=0,
\qquad
(x_{i}^-)^{\lambda(h_{i})+1}w_{\lambda,r}=0,
\qquad
hw_{\lambda,r}=\lambda(h)w_{\lambda,r},
\end{gather*}
where $i\in I$ and $h\in\mathfrak h$.
Clearly the module $\Delta(\lambda,r)$ is a~quotient of $W(\lambda,r)$ and moreover $V(\lambda,r)$ is the unique
irreducible quotient of $W(\lambda,r)$.
It is known (see~\cite{CFK} or~\cite{CPweyl}) that $W(0,r)\cong\mathbb C$ and that, if $\lambda\ne 0$, the modules
$W(\lambda,r)$ are inf\/inite-dimensional and satisfy $\wt W(\lambda,r)\subset\conv W\lambda$
and $W(\lambda,r)[s]\ne 0$ if\/f $s\ge r$, from which we see that $ W(\lambda,r)\notin \Ob \mathcal G$.
It follows that, if we set
\begin{gather*}
\nabla(\lambda,r)= W(-w_\circ\lambda,-r)^*,
\end{gather*}
then $\nabla(\lambda, r)\in\Ob {\mathcal G}_{\bdd}$ and
$\soc\nabla(\lambda,r)\cong V(\lambda, r)$.

We note that $\Delta(\lambda, r)$ (resp.
$\nabla(\lambda, r)$) is the maximal quotient of $P(\lambda, r)$ (resp.
submodule of $I(\lambda, r)$) satisfying
\begin{gather*}
[\Delta(\lambda,r):V(\mu,s)] \ne 0\implies (\mu,s)\le (\lambda,r),
\\
\big(\text{resp.}
\quad
[\nabla(\lambda,r):V(\mu,s)] \ne 0\implies (\mu,s)\le (\lambda,r) \big).
\end{gather*}
Hence these are the standard (resp.
costandard) modules in $\mathcal G$ associated to $(\lambda, r)$.

\subsection{Truncated subcategories}\label{section2.5}

In this section, we recall the def\/inition of certain Serre subcategories of
$\widehat{\mathcal G}$.

Given $\Gamma \subset \Lambda$, let $\widehat{\mathcal G}(\Gamma)$ be the full subcategory of $\widehat{\mathcal G}$
consisting of all $M$ such that
\begin{gather*}
M \in \Ob\widehat{\mathcal G},
\qquad
[M:V(\lambda,r)] \ne 0 \implies (\lambda,r)\in \Gamma.
\end{gather*}
The subcategories $\mathcal{G}(\Gamma)$ and $\mathcal G_{\bdd}(\Gamma)$ are def\/ined in the obvious way.
Observe that if $(\lambda,r)\in\Gamma$, then $V(\lambda,r) \in \widehat{\mathcal G}(\Gamma)$, and we have the following
trivial result.

\begin{Lemma}
The isomorphism classes of simple objects of $\widehat{\mathcal G}(\Gamma)$ are indexed by $\Gamma$.
\end{Lemma}

\begin{Remark}
Let $\mathcal C$ be one of the categories $\mathcal G_{\le s}$, $\mathcal G$, $\mathcal G_{\bdd}$, $\mathcal
G(\Gamma)$, $\mathcal G_{\bdd}(\Gamma)$, which are full sub\-ca\-te\-gories of $\widehat{\mathcal G}$.
Then, we have $\Ext^1_{\widehat{\mathcal G}}(M,N) = \Ext^1_{\mathcal C}(M,N)$
for all $M,N\in \mathcal C$.
\end{Remark}

\subsection{A specif\/ic truncation}\label{section2.6}

We now focus on $\Gamma$ of the form $\Gamma = P^+ \times J$, where $J$ is an interval in $\mathbb Z$ with one of the
forms $(-\infty,n]$, $[m,n]$, $[m,\infty)$ or $\mathbb Z$, where $n,m \in \mathbb Z$.
We set $a=\inf  J$ and $b=\sup  J$.
Throughout this section, we assume that $(\lambda, r)\in \Gamma$.

Let $P(\lambda,r)(\Gamma)$ be the maximal quotient of $P(\lambda,r)$ which is an object of $\mathcal G(\Gamma)$ and let
$I(\lambda,r)(\Gamma)$ be the maximal submodule of $I(\lambda,r)$ which is an object of $\mathcal G(\Gamma)$.
These are the indecomposable projective and injective modules associated to the simple module $V(\lambda,r)\in \mathcal
G(\Gamma)$.

For an object $M\in \mathcal G$, let $M^\Gamma$ be the subquotient
\begin{gather*}
M^\Gamma = \frac{M_{\ge a}}{M_{>b}},
\end{gather*}
where $M_{\ge a}= \bigoplus_{r\ge a} M[r]$, and we understand $M_{\ge a}=M$ if $a=-\infty$ and $M_{>b}=0$ if $b=\infty$.

\begin{Remark}\quad
\begin{enumerate}\itemsep=0pt
\item
If $M=\bigoplus_{s\ge p} M[s]$ for some $p \ge a$, then $M^\Gamma = \frac{M}{M_{>b}}$.
\item
If $M=\bigoplus_{s<p} M[s]$ for some $p \le b$, then $M^\Gamma = M_{\ge a}$.
\end{enumerate}
\end{Remark}

Clearly $M^\Gamma \in \mathcal G(\Gamma)$, and because morphisms are graded, this assignment def\/ines a~functor from
$\mathcal G$ to $\mathcal G(\Gamma)$.
It follows from Lemma~\ref{functor} that $\Gamma$ is exact.

If we def\/ine another subset $\Gamma'=P^+ \times \{ -J\}$, then it follows from the def\/inition of the graded duality that
if $M\in\Ob \mathcal G(\Gamma)$ then $M^*\in\Ob\mathcal G(\Gamma')$.

\begin{Lemma}
The module $P(\lambda, r)(\Gamma)=P(\lambda, r)^\Gamma$ and $I(\lambda, r)(\Gamma)=I(\lambda,r)^\Gamma$.
\end{Lemma}

We set
\begin{gather*}
\Delta(\lambda, r)(\Gamma):=\Delta(\lambda, r)^{\Gamma},
\qquad
W(\lambda, r)(\Gamma):=W(\lambda, r)^\Gamma
\qquad
\text{and}
\qquad
\nabla(\lambda, r)(\Gamma):=\nabla(\lambda, r)^{\Gamma}.
\end{gather*}
In light of the above remark, we can see that
\begin{gather*}
\Delta(\lambda, r)(\Gamma)=\frac{\Delta(\lambda, r)}{\bigoplus_{s> b}\Delta(\lambda, r)[s]},
\qquad
W(\lambda, r)(\Gamma) = \frac{W(\lambda, r)}{\bigoplus_{s> b}W(\lambda, r)[s]}\qquad \text{and}
\\
\nabla(\lambda, r)(\Gamma)=\nabla(\lambda, r)_{\ge a}.
\end{gather*}
Note that, with respect to the partial order $\le$, for each $(\lambda,r)\in \Gamma$ we have $\Delta(\lambda,r)(\Gamma)$
the maximal quotient of $P(\lambda,r)(\Gamma)$ such that
\begin{gather*}
[\Delta(\lambda,r)(\Gamma):V(\mu,s)] \ne 0\implies (\mu,s)\le (\lambda,r).
\end{gather*}
Similarly, we see that $\nabla(\lambda,r)(\Gamma)$ is the maximal submodule of $I(\lambda,r)(\Gamma)$ satisfying
\begin{gather*}
[\nabla(\lambda,r)(\Gamma):V(\mu,s)] \ne 0\implies (\mu,s)\le (\lambda,r).
\end{gather*}
These modules $\Delta(\lambda,r)(\Gamma)$ and $\nabla(\lambda,r)(\Gamma)$ are called, respectively, standard and
co-standard modules associated to $(\lambda,r)\in \Gamma$.

The following proposition summarizes the properties of $\Delta(\lambda,r)(\Gamma)$ which are necessary for this paper.
They can easily be derived from the properties of the functor $\Gamma$.

\begin{Proposition}\quad
\begin{enumerate}\itemsep=0pt
\item[$1.$]
The module $\Delta(\lambda, r)(\Gamma)$ is generated as a~$\mathfrak g[t]$-module by an element $w_{\lambda,r}$ with
relations:
\begin{gather*}
\mathfrak n^+[t]w_{\lambda,r}=0,
\qquad
(x_{i}^-)^{\lambda(h_{i})+1}w_{\lambda,r}=0,
\\
(h\otimes t^s)w_{\lambda,r}= \delta_{s,0}\lambda(h)w_{\lambda,r},
\qquad
\textbf U(\mathfrak g[t])[p]w_{\lambda,r}=0,
\qquad
\text{if}
\quad p> b-r,
\end{gather*}
for all $i\in I$, $h\in\mathfrak h$ and $s\in\mathbb Z_+$, where if $b=\infty$, then the final relation is empty
relation.
\item[$2.$]
The module $\Delta(\lambda,r)(\Gamma)$ is indecomposable and finite-dimensional and, hence, an object of $\mathcal
G_{\bdd}(\Gamma)$.
\item[$3.$]
$\dim \Delta(\lambda,r)(\Gamma)_\lambda =\dim\Delta(\lambda,r)(\Gamma)[r]_\lambda=1$.
\item[$4.$]
$\wt\Delta(\lambda,r)(\Gamma)\subset\conv W\lambda$.
\item[$5.$]
The module $V(\lambda,r)$ is the unique irreducible quotient of $\Delta(\lambda,r)(\Gamma)$.
\item[$6.$]
$\{\ch_{\gr}\Delta(\lambda,r)(\Gamma):(\lambda,r)\in \Gamma\}$ is a~linearly independent subset of
$\mathbb Z[P][u,u^{-1}]$.
\end{enumerate}
\end{Proposition}

\subsection{The truncated global Weyl modules}\label{section2.7}

Here we collect the results on $W(\lambda, r)(\Gamma)$ which we will need for this paper.
\begin{Proposition}\quad
\begin{enumerate}\itemsep=0pt
\item[$1.$]
The module $W(\lambda, r)(\Gamma)$ is generated as a~$\mathfrak g[t]$-module by an element $w_{\lambda,r}$ with
relations:
\begin{gather*}
\mathfrak n^+[t]w_{\lambda,r}=0,
\qquad
(x_{i}^-)^{\lambda(h_{i})+1}w_{\lambda,r}=0,
\qquad
h w_{\lambda,r}=\lambda(h)w_{\lambda,r},
\\
\textbf U(\mathfrak g[t])[p]w_{\lambda,r}=0,
\qquad
\text{if}
\quad
p> b-r,
\end{gather*}
where, if $b=\infty$, then the final relation is empty relation.
Here $i\in I$ and $h\in\mathfrak h$.
\item[$2.$]
The module $W(\lambda,r)(\Gamma)$ is indecomposable and an object of $\widehat{\mathcal G}(\Gamma)$.
\item[$3.$]
$\dim W(\lambda,r)(\Gamma)[r]_\lambda=1$ and $\dim W(\lambda, r)(\Gamma)_\lambda[s]\ne 0$ if and only if $r\le s \le b$.
\item[$4.$]
$\wt W(\lambda,r)(\Gamma)\subset\conv W\lambda$.
\item[$5.$]
The module $V(\lambda,r)$ is the unique irreducible quotient of $W(\lambda,r)(\Gamma)$.
\item[$6.$]
$\{ \ch_{\gr}W(\lambda,r)(\Gamma):(\lambda,r)\in \Gamma \}$ is a~linearly independent subset of
$\mathbb Z[P][u,u^{-1}]$.
\end{enumerate}
\end{Proposition}

\subsection{The costandard modules}\label{section2.8}

The following proposition summarizes the main results on $\nabla(\lambda,r)(\Gamma)$ that are needed for this paper.
All but the f\/inal result can be found by considering the properties of the functor $\Gamma$ and the paper~\cite{bc}.

\begin{Proposition}\quad
\begin{enumerate}\itemsep=0pt
\item[$1.$]
The module $\nabla(\lambda, r)(\Gamma)$ is an indecomposable object of $\mathcal G_{\bdd}(\Gamma)$.
\item[$2.$]
$\dim\nabla(\lambda,r)(\Gamma)[r]_\lambda=1$ and $\dim\nabla(\lambda,r)(\Gamma)[s]_\lambda\ne 0 \Leftrightarrow a\le
s\le r$.
\item[$3.$]
$\wt\nabla(\lambda,r)(\Gamma)\subset\conv W\lambda$.
\item[$4.$]
Any submodule of $\nabla(\lambda,r)(\Gamma)$ contains $\nabla(\lambda,r)(\Gamma)[r]_\lambda$ and the socle of
$\nabla(\lambda,r)(\Gamma)$ is the simple module $V(\lambda,r)$.
\item[$5.$]
$\{\ch_{\gr}\nabla(\lambda,r)(\Gamma):(\lambda,r)\in \Gamma \}$ is a~linearly independent subset of
$\mathbb Z[P][u,u^{-1}]]$.
\item[$6.$]
Let $\Gamma' = P^+ \times\{ - J\}$ and $(\lambda,r)\in \Gamma$.
Then $\nabla(\lambda, r)(\Gamma)\cong W(-\omega_0 \lambda, -r)(\Gamma')^*$.
\end{enumerate}
\end{Proposition}

\begin{proof} We prove the f\/inal item.
As a~vector space we have
\begin{gather*}
W(-\omega_0 \lambda, -r)(\Gamma ')\cong \bigoplus_{s=-r}^{-a}W(-\omega_0 \lambda, -r)[s].
\end{gather*}
Since $W(-\omega_0 \lambda, -r)(\Gamma ')$ is a~quotient of $W(-\omega_0 \lambda, -r)$, its dual must be a~submodule of
$\nabla(\lambda, r)$.
By the def\/inition of the graded dual, we see that, as a~vector space,
\begin{gather*}
W(-\omega_0 \lambda, -r)(\Gamma ')^* \cong \bigoplus_{s=a}^{r}\nabla(\lambda, r)[s].
\end{gather*}
Hence, as vector spaces, we see that $\nabla(\lambda, r)(\Gamma)\cong W(-\omega_0 \lambda, -r)(\Gamma')^*$.
Now, the fact that $ W(-\omega_0 \lambda, -r)(\Gamma')^*$ is a~submodule completes the proof.
\end{proof}

\section{The main theorem and some homological results}\label{section3}

\begin{Definition}
We say that $M\in\Ob\mathcal G(\Gamma)$ admits a~$\Delta(\Gamma)$ (resp.
$\nabla(\Gamma)$)-f\/iltration if there exists an increasing family of submodules $0\subset M_0\subset M_1\subset \cdots$
with $M=\bigcup_k M_k$, such that
\begin{gather*}
M_k/M_{k-1}\cong \bigoplus_{(\lambda,r)\in \Gamma}\Delta(\lambda,r)(\Gamma)^{m_k(\lambda,r)}
\qquad
\left(\text{resp.}\
\
M_k/M_{k-1}\cong \bigoplus_{(\lambda,r)\in \Gamma}\nabla(\lambda,r)(\Gamma)^{m_k(\lambda,r)}\right)
\end{gather*}
for some choice of $m_k(\lambda,r)\in \mathbb Z_+$.
We do not require $\sum\limits_{(\lambda, r)}m_k(\lambda, r)<\infty$.
If $M_k=M$ for some $k\ge 0$, then we say that $M$ admits a~f\/inite $\Delta(\Gamma)$ (resp.\
$\nabla(\Gamma)$)-f\/iltration.
Because our modules have f\/inite-dimensional graded components, we can conclude that the multiplicity of a~f\/ixed
$\Delta(\lambda, r)(\Gamma)$ (resp.\
$\nabla(\lambda, r)(\Gamma))$ in a~$\Delta(\Gamma)$-f\/iltration (resp.\
$\nabla(\Gamma)$-f\/iltration) must be f\/inite, and we denote this multiplicity by $[M:\Delta(\lambda, r)(\Gamma)]$ (resp.
$[M:\nabla(\lambda, r)(\Gamma)]$).
Finally, we say that $M\in\Ob\mathcal G(\Gamma)$ is tilting if $M$ has both a~$\Delta(\Gamma)$ and
a~$\nabla(\Gamma)$-f\/iltration.
\end{Definition}

The main goal of this paper is to understand tilting modules in $\mathcal G_{\bdd}(\Gamma)$.
(The case where $J=\mathbb Z$ was studied in~\cite{bc}.) In the case of algebraic groups (see~\cite{Donkin,Mathieu})
a~crucial necessary result is to give a~cohomological characterization of modules admitting
a~$\nabla(\Gamma)$-f\/iltration.
The analogous result in our situation is to prove the following statement:

\textit{An object $M$ of $\mathcal G_{\bdd}(\Gamma)$ admits a~$\nabla(\Gamma)$-filtration if and only if
$\Ext^1_{{\widehat{\mathcal G}}}((\Delta(\lambda,r)(\Gamma), M)=0$ for all $(\lambda,r)\in \Gamma$.}

It is not hard to see that the forward implication is true.
The converse statement however requires one to prove that any object of $\mathcal G_{\bdd}(\Gamma)$ can be
embedded in a~module which admits a~$\nabla(\Gamma)$-f\/iltration.
This in turn requires Theorem~\ref{bgg}.
Summarizing, the f\/irst main result that we shall prove in this paper is:

\begin{Proposition}
\label{extnablaconn}
Let $M\in\Ob\mathcal G_{\bdd}(\Gamma)$.
Then the following are equivalent:{\samepage
\begin{enumerate}\itemsep=0pt
\item[$1.$]
The module $M$ admits a~$\nabla(\Gamma)$-filtration.
\item[$2.$]
$M$ satisfies $\Ext^1_{\widehat{\mathcal G}}(\Delta(\lambda,r)(\Gamma), M)=0$
for all
$(\lambda,r)\in \Gamma$.
\end{enumerate}}
\end{Proposition}

The second main result that we shall prove in this paper is the following:
\begin{Theorem}\quad
\begin{enumerate}\itemsep=0pt
\item[$1.$]
Given $(\lambda, r)\in\Gamma$, there exists an indecomposable module $T(\lambda,r)(\Gamma)\in\Ob\mathcal
G_{\bdd}(\Gamma)$ which admits a~$\Delta(\Gamma)$-filtration and a~$\nabla(\Gamma)$-filtration.
Further,
\begin{gather*}
T(\lambda,r)(\Gamma)[s]_\lambda=0
\qquad
\text{if}
\quad
s>r,
\quad
T(\lambda,r)(\Gamma)[r]_\lambda=1,
\quad
\wt T(\lambda,r)(\Gamma)\subset\conv W\lambda,
\end{gather*}
and $T(\lambda,r)(\Gamma)\cong T(\mu,s)(\Gamma)$ if and only if $(\lambda,r)=(\mu,s)$.
\item[$2.$]
Moreover, any indecomposable tilting module in $\mathcal G_{\bdd}(\Gamma)$ is isomorphic to
$T(\lambda,r)(\Gamma)$ for some $(\lambda,r)\in \Gamma$, and any tilting module in $\mathcal
G_{\bdd}(\Gamma)$ is isomorphic to a~direct sum of indecomposable tilting modules.
\end{enumerate}
\end{Theorem}

\section{Proof of Proposition~\ref{extnablaconn}}\label{section4}

\subsection{Initial homological results}\label{section4.1}

 We begin by proving the implication (1) $\implies$ (2) from
Proposition~\ref{extnablaconn}.
In order to do this, we f\/irst establish some homological properties of the standard and costandard modules which will be
used throughout the paper.

\begin{Proposition}
\label{homresults}
Let $\lambda, \mu \in P^+$.
Then we have the following
\begin{enumerate}\itemsep=0pt
\item[$1.$]
$\Ext^1_{\widehat{\mathcal G}}(W(\lambda,r)(\Gamma),
W(\mu,s)(\Gamma))=0=\Ext^1_{\widehat{\mathcal G}}(\nabla(\mu,s)(\Gamma), \nabla(\lambda,r)(\Gamma)) $ for
all $s,r\in \mathbb Z$ if $\lambda \not < \mu$.
\item[$2.$]
$\Ext^1_{\widehat{\mathcal G}}(\Delta(\lambda, r)(\Gamma), \nabla(\mu,s)(\Gamma))=0$.
\item[$3.$]
If $\lambda \not \le \mu$ then $\Ext^1_{\widehat{\mathcal G}}(\Delta(\lambda, r)(\Gamma),
\Delta(\mu,s)(\Gamma))=0$ for all $s,r\in \mathbb Z$.
\item[$4.$]
If $s\ge r$, then $\Ext^1_{\widehat{\mathcal G}}(\Delta(\lambda, s)(\Gamma), \Delta(\lambda,r)(\Gamma))=0$.
\end{enumerate}
\end{Proposition}

\begin{proof} For part (1), suppose that we have a~sequence $0\rightarrow W(\mu,s)(\Gamma) \rightarrow X \rightarrow
W(\lambda,r)(\Gamma) \rightarrow 0$, and let $x\in X_\lambda[r]$ be a~pre-image of $w_{\lambda,r}$.
It is clear from the hypothesis on $\mu$ that $\mathfrak n^+[t].x=0$.
If $b< \infty$ and $p>b-r$, then $\textbf U(\mathfrak g[t])[p].x=0$ by grade considerations.
Note that, since $\dim X[r]<\infty$, we have $\dim \textbf U(\mathfrak g)x<\infty$.
It follows from the f\/inite-dimensional representation theory that $(x_i^-)^{\lambda(h_i)+1}x=0$, for all $i\in I$, and
so the sequence splits.
The proof for $\nabla(\Gamma)$ is similar and omitted.

For part (2), suppose we have a~sequence $0\rightarrow \nabla(\mu,s)(\Gamma) \rightarrow X \rightarrow
\Delta(\lambda,r)(\Gamma) \rightarrow 0$ and $\mu\not\ge \lambda$.
Then $\dim X_\lambda[r]=1$, and if $x\in X_\lambda[r]$ is a~pre-image of $w_{\lambda,r}$ then $x$ satisf\/ies the def\/ining
relations of $w_{\lambda,r}$ and the sequence splits.
If $\mu \ge \lambda$, then by taking duals we get a~sequence $0\rightarrow \Delta(\lambda,r)(\Gamma)^*\rightarrow Y
\rightarrow W(-\omega_0\mu,-s)(\Gamma')\rightarrow 0$.
Again, if $y\in Y_{-\omega_0 \mu}[-s]$ is a~pre-image of $w_{-\omega_0 \mu, -s}$, we see that $y$ satisf\/ies the def\/ining
relations of $w_{-\omega_0 \mu,-s}\in W(-\omega_0\mu,-s)(\Gamma')$, and the sequence splits.

The proofs for parts (3) and (4) are similar to that for part (1) and are omitted.
\end{proof}

The proof of the following lemma is standard (see, for example,~\cite{bc}).

\begin{Lemma}
\label{reorder}
Suppose that $M\in \Ob \mathcal G_{\bdd}(\Gamma)$ admits a~$($possibly infinite$)$
$\nabla(\Gamma)$-filtration.
Then~$M$ admits a~finite $\nabla(\Gamma)$-filtration
\begin{gather*}
0\subset M_1\subset M_2 \subset\dots\subset M_k=M
\qquad
\text{with}
\quad
M_s / M_{s-1} \cong \bigoplus_{r\in \mathbb Z} \nabla(\lambda_s, r)(\Gamma)^{[M:\nabla(\lambda_s,r)(\Gamma)]},
\end{gather*}
where $\lambda_i>\lambda_j$ implies $i>j$.
In particular if $\mu$ is maximal such that $M_\mu \ne 0$, then there exists $s\in\mathbb Z$ and a~surjective map
$M\rightarrow \nabla(\mu,s)(\Gamma)$ such that the kernel has a~$\nabla(\Gamma)$-filtration.
\end{Lemma}

We can now prove the implication (1) $\implies$ (2) from Proposition~\ref{extnablaconn}.

\begin{Corollary}
If $M\in \Ob\mathcal G_{\bdd}(\Gamma)$ admits a~$\nabla(\Gamma)$-filtration, then
$\Ext^1_{\mathcal G(\Gamma)}(\Delta(\lambda,r)(\Gamma), M)=0$, for all $(\lambda,r)\in \Gamma$.
\end{Corollary}

\begin{proof} Let $0\subset M_1 \subset M_2 \subset\dots\subset M_k=M$ be a~f\/inite $\nabla(\Gamma)$-f\/iltration as in
Lemma~\ref{reorder}.
Then, $M_s / M_{s-1}$ is a~(possibly inf\/inite) direct sum of $\nabla(\lambda_s, r)(\Gamma)$, and
\begin{gather*}
\Ext^1_{\mathcal G(\Gamma)}(\Delta(\lambda, r)(\Gamma), M_s / M_{s-1})
\cong
\Ext^1_{\mathcal G(\Gamma)}\left(\Delta(\lambda, r)(\Gamma), \bigoplus_{r\in \mathbb Z} \nabla(\lambda_s,r)(\Gamma)^{m_s(\lambda_s,r)}\right)
\\
\hphantom{\Ext^1_{\mathcal G(\Gamma)}(\Delta(\lambda, r)(\Gamma), M_s / M_{s-1})}{}
\cong \prod \Ext^1_{\mathcal G(\Gamma)}(\Delta(\lambda, r)(\Gamma), \nabla(\lambda_s, r)(\Gamma)) =0,
\end{gather*}
by Proposition~\ref{homresults} (2).
The result follows by induction on $k$, the length of the f\/iltration.
\end{proof}

\subsection{Towards understanding extensions between the standard\\ and costandard modules}\label{section4.2}

\begin{Proposition}
Suppose that $N\in\Ob\mathcal G(\Gamma)$ is such that $\Ext^1_{\mathcal
G(\Gamma)}(\Delta(\lambda,r)(\Gamma), N)=0$ for all $(\lambda,r)\in \Gamma$.
If $M\in\Ob\mathcal G(\Gamma)$ has a~$\Delta(\Gamma)$-filtration then $\Ext^1_{\mathcal G(\Gamma)}(M, N)=0$.
\end{Proposition}

\begin{proof} Consider a~short exact sequence $0\to N\to U\to M\to 0$.
Suppose that $M_k\subset M_{k+1}$ is a~part of the $\Delta(\Gamma)$-f\/iltration of $M$ and assume that
\begin{gather*}
M_{k+1}/M_k\cong\bigoplus_{(\mu,s)\in\Lambda}\Delta(\mu,s)(\Gamma)^{m_k(\mu,s)}.
\end{gather*}
By assumption we have $\Ext^1_{\mathcal G(\Gamma)}(M_{k+1}/M_k, N)=0$.
Let $U_k\subset U$ be the pre-image of $M_k$, which contains $N$ because $0\in M_k$.
Note that $U_{k+1}/U_k\cong M_{k+1}/M_k$.
Now, consider the short exact sequence $0\to N\to U_k\to M_k\to 0$.
This sequence def\/ines an element of $\Ext^1_{\mathcal G(\Gamma)}(M_k, N)$.
Since~$M_k$ has a~f\/inite $\Delta(\Gamma)$-f\/iltration it follows that $\Ext^1_{\mathcal G(\Gamma)}(M_k,
N)=0$.
Hence the sequence splits and we have a~retraction $\varphi_k: U_k\to N$.
We want to prove that $\varphi_{k+1}:U_{k+1}\to N$ can be chosen to extend $\varphi_k$.
For this, applying $\Hom_{\mathcal G(\Gamma)}(-,N)$ to $0\to U_k\to U_{k+1}\to U_{k+1}/U_k\to 0$, we get
$\Hom_{\mathcal G(\Gamma)}(U_{k+1}, N)\to \Hom_{\mathcal G(\Gamma)}(U_k, N)\to 0$, which
shows that we can choose $\varphi_{k+1}$ to lift~$\varphi_k$.
Now def\/ining $\varphi: U\to N$ by $\varphi(u)=\varphi_k(u)$, for all $u\in U_k$, we have the desired splitting of the
original short exact sequence.
\end{proof}

Together with Proposition~\ref{homresults} and taking $N=\nabla(\lambda,r)(\Gamma)$ in the proposition above, we now
have:

\begin{Corollary}
\label{extMnabla}
Suppose $M\!\in\!\Ob\mathcal G(\Gamma)$ admits a~$\Delta(\Gamma)$-filtration.
Then, $\Ext^1_{\widehat{\mathcal G}}(M,\nabla(\lambda,r)(\Gamma))\!=\!0$, for all $(\lambda,r)\in\Gamma$.
\end{Corollary}

\subsection{A natural embedding}\label{section4.3}

In this section we show that every $M\in \Ob \mathcal G_{\bdd}(\Gamma)$ embeds into an
injective module $I(M)\in \Ob \mathcal G(\Gamma)$.
Let $\soc M\subset M$ be the maximal semi-simple submodule of $M$.

\begin{Lemma}\label{embeda}
Let $M\in \Ob \mathcal G_{\bdd}(\Gamma)$.{\samepage
\begin{enumerate}\itemsep=0pt
\item[$1.$]
If $M\ne 0$, then $\soc M \ne 0$.
\item[$2.$]
Suppose $\soc M=\bigoplus V(\lambda, r)^{m_{\lambda, r}}$.
Then $M \hookrightarrow \bigoplus I(\lambda,r)(\Gamma)^{m_{\lambda,r}}$.
\end{enumerate}}
\end{Lemma}

\begin{proof} For the f\/irst part, let $M\in \Ob\mathcal G_{\bdd}$, suppose $M\ne 0$, and let
$s\in\mathbb Z$ be minimal such that $M\in \Ob \mathcal G_{\le s}$.
Then $M[s]\ne 0$ and $M[s]\subset \soc M$.

For the second, let $\soc M=\bigoplus_{(\lambda, r)\in \Lambda} V(\lambda, r)^{m_{\lambda, r}}$ from which
we get a~natural injection $\soc M \stackrel{ \iota}\hookrightarrow \bigoplus
I(\lambda,r)^{m_{\lambda,r}}$.
By injectivity, we have a~morphism $ M \stackrel{f}\rightarrow \bigoplus I(\lambda,r)^{m_{\lambda,r}}$, through which~$\iota$ factors.
In fact, we can show that $f$ is an injection.
If not, $\soc \ker f \ne 0$.
On the other hand, $\soc \ker f \subset \soc M$, and $f$ is injective on $\soc
M$.
If we assume that $M\in \Ob \mathcal G_{\bdd}(\Gamma)$, then it is easy to conclude that $\operatorname{im}
f \subset \bigoplus I(\lambda,r)(\Gamma)^{m_{\lambda,r}}$, completing the proof.
\end{proof}

\subsection{o-canonical f\/iltration}\label{section4.4}

In this section we shall establish a~f\/inite f\/iltration on modules $M\in \Ob \mathcal
G_{\bdd}(\Gamma)$ where the successive quotients embed into direct sums of $\nabla(\Gamma)$.
We then use the f\/iltration to establish lower and upper bounds on the graded character of $M$.
We use the character bounds to prove Proposition~\ref{extnablaconn}.

Now f\/ix an ordering of $P^+=\{\lambda_0,\lambda_1,\dots\}$ such that $\lambda_r>\lambda_s$ implies that $r>s$.
For $M\in\Ob \mathcal G(\Gamma)$ we set $M_s \subset M$ as the maximal submodule whose weights lie in \mbox{$\{
\conv W\lambda_r \,|\, r\le s\}$}.
Evidently $M_{s-1}\subset M_s$.
We call this the o-canonical f\/iltration, because it depends on the order.
This is a~f\/inite f\/iltration for $M\in\Ob \mathcal G_{\bdd}(\Gamma)$ and we set $k(M)$ to be
minimal such that $M=M_{k(M)}$.
Clearly
\begin{gather}
M_{s-1}\subset M_s,
\qquad
M=\bigcup_{s=0}^{k(M)}M_s,
\qquad
\text{and}
\nonumber
\\
\label{limit}
\Hom_{\mathcal G}(V(\lambda, r), M_s/M_{s-1})\ne 0\implies \lambda=\lambda_s.
\end{gather}
It follows from Lemma~\ref{embeda} and~\eqref{limit} that the quotient $M_s/M_{s-1}$ embeds into a~module of the form
$\bigoplus I(\lambda_s, r)(\Gamma)^{m_{s,r}}$, where $m_{s,r}=\dim \Hom_{\mathcal G}(V(\lambda_s,r),
M_s/M_{s-1})$.
Since the weights of $M_s / M_{s-1}$ are bounded above by $\lambda_s$, they embed into the maximal submodule of this
direct sum, whose weights are bounded above by $\lambda_s$.
Hence we have $M_s / M_{s-1}$ embedding into a~direct sum of modules of the form $\nabla(\lambda_s, r)$, with $r\in J$.
This gives,
\begin{gather*}
\ch_{\gr}M=\sum\limits_{s\ge 0}\ch_{\gr}M_{s}/M_{s-1}\le \sum\limits_{s\ge
0}\sum\limits_{r\in J}\dim\Hom_{\mathcal G}(V(\lambda_{s},r),
M_s/M_{s-1})\ch_{\gr}\nabla(\lambda_{s},r)(\Gamma),
\end{gather*}
i.e.,
\begin{gather*}
[M: V(\mu,\ell)]\le \sum\limits_{s\ge 0}\sum\limits_{r\in J}\dim\Hom_{\mathcal G}(V(\lambda_{s},r),
M_s/M_{s-1})[\nabla(\lambda_{s},r)(\Gamma): V(\mu,\ell)],
\end{gather*}
for all $(\mu,\ell)\in\Lambda$.
We claim that this is equivalent to
\begin{gather*}
\ch_{\gr} M = \sum \ch_{\gr} M_s / M_{s-1} \le \sum\limits_{s\ge
0}\sum\limits_{r\in J} \dim \Hom_{\mathcal G} (\Delta(\lambda_s,r)(\Gamma),
M)\ch_{\gr} \nabla(\lambda_s,r)(\Gamma).
\end{gather*}

The claim follows from Lemma~\ref{aaa} below, and~\cite[\S~3.5]{bc}, which states that
\begin{gather*}
\Hom_{\mathcal G} (\Delta(\lambda_s,r),M) \cong \Hom_{\mathcal G} \left(V(\lambda_s,r), M_s
/ M_{s-1} \right).
\end{gather*}

\begin{Lemma}
\label{aaa}
Let $M\in \Ob\mathcal G(\Gamma)$ and $(\lambda,r)\in \Gamma$.
We have
\begin{gather*}
\Hom_{\widehat{\mathcal G}}(\Delta(\lambda,r)(\Gamma),M) \cong \Hom_{\widehat{\mathcal
G}}(\Delta(\lambda,r),M)
\end{gather*}
{and}
\begin{gather*}
\Ext^1_{\widehat{\mathcal G}}(\Delta(\lambda,r)(\Gamma),M) \cong \Ext^1_{\widehat{\mathcal
G}}(\Delta(\lambda,r),M).
\end{gather*}
\end{Lemma}

\begin{proof} As $M[s]=0$ for all $s>b$, we must have $\Hom_{\widehat{\mathcal G}}(\bigoplus_{s>b}
\Delta(\lambda,r)[s],M)=0$.
Similarly, if $0\to M \to X \to \bigoplus_{s>b} \Delta(\lambda,r)[s] \to 0$ is exact, by using again that $M[s]=0$ for
$s>b$, we have $\dim X[n] = \dim (\bigoplus_{s>b}\Delta(\lambda,r)[s])[n]$ for any $n>b$ and so we have an injective map
$\iota: \bigoplus_{s>b}\Delta(\lambda,r)[s] \to X$ which splits the sequence.
Thus, $\Ext^1_{\widehat{\mathcal G}}(\bigoplus_{s>b} \Delta(\lambda,r)[s],M)=0$.
Now the other statements are easily deduced.
\end{proof}

Therefore, we get
\begin{gather}
\label{lowbd}
\ch_{\gr} M \le \sum\limits_{s\ge 0} \sum\limits_{r\in J} \dim \Hom_{\mathcal
G(\Gamma)}(\Delta(\lambda_s,r)(\Gamma),M) \ch_{\gr} \nabla(\lambda_s,r)(\Gamma)
\end{gather}
and the equality holds if, and only if, the o-canonical f\/iltration is $\nabla(\Gamma)$-f\/iltration.

\subsection{A homological characterization of costandard modules}\label{section4.5}

Note that even if $N\notin \Ob\mathcal
G_{\bdd}$, we can still def\/ine the submodules $N_s$, and by def\/inition $N_s\in \Ob\mathcal
G_{\bdd}$.
The following result, or more precisely a~dual statement about projective modules and global Weyl f\/iltrations, was
proved for $\mathfrak g= \mathfrak{sl}_2$ in~\cite{bcm}, for $\mathfrak g =\mathfrak {sl}_{n+1}$ in~\cite{bbcklk} and
for general $\mathfrak g$ in~\cite{ci}.
In particular, we note that the argument in \cite[Section~5.5]{bbcklk} works in general.

\begin{Theorem}
\label{bgg}
For all $(\lambda, r) \in \Lambda$ and for all $p\in \mathbb Z_+$ the o-canonical filtration on $I(\lambda,r)_p$ is
a~$\nabla$-filtration.
Moreover, for all $(\mu,s)$ we have $[I(\lambda,r)_p:\nabla(\mu,s)]=[\Delta(\mu,r):V(\lambda,s)]$.
\end{Theorem}

We combine this with equation~\eqref{lowbd} and the linear independence of the graded characters of the $\nabla(\lambda,
r)(\Gamma)$ to see that $[I(\lambda,r)_p:\nabla(\mu,s)]=\dim \Hom_{\mathcal
G(\Gamma)}(\Delta(\lambda_s,r),I(\lambda, r)_p)$.

As a~consequence of the theorem and the exactness of the functor $\Gamma$, we conclude that $I(\lambda,r)_p(\Gamma)$ has
a~$\nabla(\Gamma)$-f\/iltration.
It is easy to see that $I(\lambda,r)_p(\Gamma)\in \mathcal G_{\bdd}(\Gamma)$.

For $M\in \mathcal G_{\bdd}(\Gamma)$, let $p$ be minimal such that $M_p=M$.
Then it is clear that we can ref\/ine the embedding from Lemma~\ref{embeda} to $M \hookrightarrow \bigoplus_{}
I(\lambda,r)_p(\Gamma)$.
We can conclude the following:

\begin{Corollary}
\label{embed}
For all $M\in \Ob\mathcal G_{\bdd}(\Gamma)$, we have $M\subset I(M) \in
\Ob\mathcal G_{\bdd}(\Gamma)$, where the o-canonical filtration of $I(M)$ is
a~$\nabla(\Gamma)$-filtration.
\end{Corollary}

We now complete the proof of Proposition~\ref{extnablaconn} following the argument in~\cite{bc}.
Let $M \in \mathcal G_{\bdd}(\Gamma)$ and assume that $\Ext^1_{\widehat{\mathcal
G}}(\Delta(\lambda,r)(\Gamma),M)=0$.
Let $I(M)\in \mathcal G_{\bdd}(\Gamma)$ be as in Corollary~\ref{embed}, and consider the short exact
sequence $0 \to M \to I(M) \to Q \to 0$.
where $Q \in \Ob \mathcal G_{\bdd} (\Gamma)$.
The assumption on the module $M$ implies that, if we apply the functor $\Hom_{\mathcal
G(\Gamma)}(\Delta(\lambda,r)(\Gamma), -)$, we get the short exact sequence
\begin{gather}
0
\rightarrow
\Hom_{\mathcal G(\Gamma)}(\Delta(\lambda,r)(\Gamma), M)
\rightarrow
\Hom_{\mathcal G(\Gamma)}(\Delta(\lambda,r)(\Gamma), I(M))  \nonumber\\
\hphantom{0}{}
\rightarrow
\Hom_{\mathcal G(\Gamma)}(\Delta(\lambda,r)(\Gamma), Q)
\rightarrow
0.\label{equal1}
\end{gather}
Since the o-canonical f\/iltration of $I(M)$ is a~$\nabla(\Gamma)$-f\/iltration, we can conclude that~\eqref{lowbd} is an
equality for $I(M)$.
We get that
\begin{gather*}
\ch_{\gr}M
= \ch_{\gr} I(M) - \ch_{\gr} Q
\ge \sum\limits_{\substack{r\in J\\s\ge 0}}
(\dim \Hom_{\mathcal G(\Gamma)}(\Delta(\lambda_s,r)(\Gamma),I(M))\\
\phantom{\ch_{\gr}M= }
 - \dim
\Hom_{\mathcal G(\Gamma)}(\Delta(\lambda_s,r)(\Gamma),Q))
\ch_{\gr}\nabla(\lambda_s,r)
\\
\phantom{\ch_{\gr}M }
= \sum\limits_{\substack{r\in J\\s\ge 0} } \dim \Hom_{\mathcal G(\Gamma)}(\Delta(\lambda_s, r)(\Gamma),M)
\ch_{\gr}\nabla(\lambda_s,r),
\end{gather*}
where the f\/inal equality is from the exactness of~\eqref{equal1}.
We now get that the character bound in~\eqref{lowbd} is an equality for $M$, and, hence, that $M$ has
a~$\nabla(\Gamma)$-f\/iltration.

Finally, we can prove the following.
\begin{Proposition}\label{26}
The following are equivalent for a~module $M\in \Ob \mathcal G_{\bdd}(\Gamma)$
\begin{enumerate}\itemsep=0pt
\item[$1.$]
For all $(\lambda, r)\in \Gamma$, we have $\Ext^1_{\widehat{\mathcal G}}(\Delta(\lambda,r)(\Gamma),M)=0$.
\item[$2.$]
$M$ admits a~$\nabla(\Gamma)$-filtration.
\item[$3.$]
The o-canonical filtration on $M$ is a~$\nabla(\Gamma)$-filtration.
\end{enumerate}
\end{Proposition}

\begin{proof} The equivalence of (1) and (2) is precisely the statement of Proposition~\ref{extnablaconn}.
Clearly $(3)$ implies $(2)$, so it is enough to show that~$(1)$ implies~$(3)$.
Assuming $(1)$, we have shown that the character bound in~\eqref{lowbd} is an equality, which is true if and only if the
o-canonical f\/iltration is a~$\nabla(\Gamma)$-f\/iltration.
\end{proof}

\subsection{Extensions between simple modules} \label{section4.6}

Our f\/inal result before constructing the tilting modules $T(\lambda,
r)(\Gamma)$ shows that the space of extensions between standard modules is always f\/inite-dimensional.
The proof is analogous to the proof in~\cite{bc}.

\begin{Proposition}
\label{fdext}
For all $(\lambda, r),(\mu,s) \in \Gamma$ we have $\dim \Ext^1_{\widehat{\mathcal G}}(\Delta(\lambda,
r)(\Gamma), \Delta(\mu,s)(\Gamma))< \infty$.
\end{Proposition}
\begin{proof}
Consider the short exact sequence $0\rightarrow K \rightarrow P(\lambda, r)\rightarrow \Delta(\lambda, r)(\Gamma)
\rightarrow 0 $ and apply the functor $Hom_{\widehat{\mathcal G}}(-,\Delta(\mu,s)(\Gamma))$.
Since $P(\lambda, r)$ is projective, we see that the result follows if $\dim \Hom_{\widehat{\mathcal
G}}(K, \Delta(\mu,s)(\Gamma))< \infty$.
Let $\ell$ be such that $\Delta(\mu,s)(\Gamma)[p]=0$ for all $p>\ell$.
Then $\Hom_{\widehat{\mathcal G}}(K_{>\ell}, \Delta(\mu,s)(\Gamma))=0$, and, hence, we have an injection
$\Hom_{\widehat{\mathcal G}}(K, \Delta(\mu,s)(\Gamma)) \hookrightarrow
\Hom_{\widehat{\mathcal G}}(\frac{K}{K_{>\ell}}, \Delta(\mu,s)(\Gamma))$.
The proposition follows because $\frac{K}{K_{>\ell}}$ is f\/inite-dimensional.
\end{proof}

\section{Construction of tilting modules}\label{section5}

\subsection{Def\/ining a~subset which can be appropriately enumerated} \label{section5.1}

In this section we construct a~family of
indecomposable modules in the category $\mathcal G_{\bdd}(\Gamma)$, denoted by
$\{T(\lambda,r)(\Gamma):(\lambda,r)\in\Gamma\}$, each of which admits a~$\Delta(\Gamma)$-f\/iltration and satisf\/ies
\begin{gather*}
\Ext^1_{\widehat{\mathcal G}}(\Delta(\mu,s)(\Gamma), T(\lambda,r)(\Gamma))=0,
\qquad
(\mu,s)\in\Gamma.
\end{gather*}
It follows that the modules $T(\lambda,r)(\Gamma)$ are tilting and we prove that any tilting module in $\mathcal
G_{\bdd}(\Gamma)$ is a~direct sum of copies of $T(\lambda,r)(\Gamma)$, $(\lambda,r)\in\Gamma$.
The construction is a~generalization of the one from~\cite{bc}, and the ideas are similar to the ones given
in~\cite{Mathieu}.
One of the f\/irst dif\/f\/iculties we encounter when trying to construct $T(\lambda,r)(\Gamma)$, using the algorithm given
in~\cite{Mathieu}, is to f\/ind a~suitable subset (depending on $(\lambda,r)$) of $\Gamma$ which can be appropriately
enumerated.
Hence we assume the following result, whose proof we postpone to Section~\ref{section5.6}.

\begin{Proposition}\label{subset}
Fix $(\lambda, r)\in \Gamma$ and assume that under the enumeration we have $\lambda=\lambda_k$.
Then there exists a~subset $\mathcal S(\lambda,r) \subset \Gamma$ such that
\begin{enumerate}\itemsep=0pt
\item[$1)$]
$(\lambda, r)\in \mathcal S(\lambda, r)$;
\item[$2)$]
there exists $r_i$ for each $i\le k$ such that $r_i\ge r$, $r_k=r$, and
\begin{gather*}
\mathcal S(\lambda, r) = \{ (\lambda_i,s) \,|\, i\le k,
s\le r_i \};
\end{gather*}

\item[$3)$]
$\Ext^1_{\widehat{\mathcal G}}(\Delta(\mu',s')(\Gamma), \Delta(\mu,s)(\Gamma))=0$ for all $(\mu,s) \in
\mathcal S(\lambda,r)$ and $(\mu',s')\notin \mathcal S(\lambda,r)$.
\end{enumerate}
Furthermore, there exists an injection $\eta: \mathcal S(\lambda,r) \to \mathbb Z_{\ge 0}$ such that for $(\mu_i,p_i) =
\eta^{-1}(i)$ we have $\Ext^1_{\widehat{\mathcal
G}}(\Delta(\mu_i,p_i)(\Gamma),\Delta(\mu_j,p_j)(\Gamma))=0$ if $i< j$, and $\Delta(\mu_j,p_j)(\Gamma)_{\mu_i}[p_i]=0$
for $i<j$.
\end{Proposition}

Without loss of generality we may assume that $\eta(\lambda, r)=0$ and the image of $\eta$ is an interval.
We need the following elementary result.

\begin{Lemma}\label{elem}
If $M,N\in\operatorname{Ob}\widehat{\mathcal G}$ are such that
$0<\dim\operatorname{Ext}^1_{\widehat{\mathcal G}}(M,N)<\infty$ and $\operatorname{Ext}^1_{\widehat{\mathcal G}}(M,M)=0$. Then, there exists $U\in\operatorname{Ob}\widehat{\mathcal G}$, $d\in\mathbb Z_+$ and a non-split exact sequence
$0\to N\to U\to M^{d}\to 0$ so that $\operatorname{Ext}^1_{\widehat{\mathcal G}}(M,U)=0$.
\end{Lemma}

\begin{proof} The proof follows from induction on $\dim\Ext^1_{\widehat{\mathcal G}}(M,N)$.
The base case is obvious.
For the inductive step, chose any non-split sequence $0\to N\to U'\to M\to 0$.
Apply the functor $\Hom_{\widehat{\mathcal G}}(M,-)$ to the sequence, and note that the image of map
$1_M\in \Hom_{\widehat{\mathcal G}}(M,M)$ is in the kernel of the surjection from
$\Ext^1_{\widehat{\mathcal G}}(M,N) \rightarrow \Ext^1_{\widehat{\mathcal G}}(M,U')$.
It follows that $\dim\Ext^1_{\widehat{\mathcal G}}(M,U') < \dim\Ext^1_{\widehat{\mathcal
G}}(M,N)$.
By the induction hypothesis, we now have a~module $U$ and a~non-split sequence $0 \to U' \to U \to M^{d-1} \to 0$.
Now, considering the sequence $0 \to U'/N \to U/N \to M^{d-1} \to 0$, and again using that
$\Ext^1_{\widehat{\mathcal G}}(M,M)=0$, we get a~non-split sequence $0 \to N \to U \to M^{d} \to 0$.
\end{proof}

\subsection{Constructing tilting modules}\label{section5.2}

 We now use $\eta$ to construct an inf\/inite family of f\/inite-dimensional
modules $M_i$, whose direct limit will be $T(\lambda, r)(\Gamma)$.
We note that the construction, at this point, will seem to be dependent on the ordering of $P^+$ we have chosen, and on
the set $\mathcal S(\lambda, r)$ and $\eta$.
We prove independence at the end of this section.

Set $M_0=\Delta(\mu_0, p_0)(\Gamma)$.
If $\Ext^1_{\widehat{\mathcal G}}(\Delta(\mu_1,p_1)(\Gamma), \Delta(\mu_0,p_0)(\Gamma))=0$, then set
$M_1=M_0$.
If not, then since $\dim(\Ext^1_{\widehat{\mathcal G}}(\Delta(\mu_1,p_1)(\Gamma),
\Delta(\mu_0,p_0)(\Gamma)))<\infty$ by Proposition~\ref{fdext}, Lemma~\ref{elem} gives us an object $ M_1' \in
\Ob \widehat{\mathcal G}(\Gamma)$ and a~non-split short exact sequence
\begin{gather*}
0\to M_0\to M_1' \to {\Delta(\mu_1,p_1)(\Gamma)^{d_1'}} \to 0
\end{gather*}
with $\Ext^1_{\widehat{\mathcal G}}(\Delta(\mu_1,p_1)(\Gamma),M_1')=0$.

Let $M_1\subseteq M_1'$ be an indecomposable summand containing $(M_1')_{\mu_0}[p_0]$.
By Proposition~\ref{subset} we see that $(M_1')_{\mu_0}[p_0]=(M_0)_{\mu_0}[p_0]$.
Then we have $M_0 \stackrel{\iota_0}{\hookrightarrow} M_1$.
Now, suppose that $M_1 \ne M_1'$.
Then, since $M_1'$ is generated by $(M_1')_{\mu_0}[p_0]$ and $(M_1')_{\mu_1}[p_1]$, we must have $(M_1)_{\mu_1}[p_1] \ne
(M_1')_{\mu_1}[p_1]$.
We must have $\dim (M_1')_{\mu_1}[p_1] - \dim (M_1)_{\mu_1}[p_1]$ linearly independent vectors in $ (M_1')_{\mu_1}[p_1]$
which do not have a~pre-image in $M_0$, and each one must then generate a~copy of $\Delta(\mu_1,p_1)(\Gamma)$.
So, $M_1'=M_1 \oplus \Delta(\mu_1,p_1)(\Gamma)^{ d}$ for some $d$ and we obtain the sequence
\begin{gather}
\label{eq422}
0\to M_0\to M_1 \to {\Delta(\mu_1,p_1)(\Gamma)^{ d_1}} \to 0.
\end{gather}
By applying the $\Hom(\Delta(\mu,p)(\Gamma),-) $ in the sequence~\eqref{eq422}, we get
\begin{gather*}
\cdots\rightarrow \Ext^1_{\widehat{\mathcal G}(\Gamma)}(\Delta(\mu,p)(\Gamma),M_0 )\rightarrow
\Ext^1_{\widehat{\mathcal G}(\Gamma)}(\Delta(\mu,p)(\Gamma),M_1 )
\\
\hphantom{\cdots}{}
\rightarrow\Ext^1_{\widehat{\mathcal G}(\Gamma)}(\Delta(\mu,p)(\Gamma),\Delta(\mu_1, p_1)(\Gamma)^{ d_1})\rightarrow \cdots.
\end{gather*}

If $(\mu,p)\notin \mathcal S(\lambda, r)$ or $(\mu, p)=(\mu_0, p_0)$ by Proposition~\ref{subset} we get
\begin{gather*}
\Ext^1_{\widehat{\mathcal G}(\Gamma)}(\Delta(\mu,p)(\Gamma),M_0 ) =0,
\qquad
\Ext^1_{\widehat{\mathcal G}(\Gamma)}(\Delta(\mu,p)(\Gamma),\Delta(\mu_1, p_1)(\Gamma))=0
\end{gather*}
and we see that the middle term is also trivial, i.e.
\begin{gather*}
\Ext^1_{\widehat{\mathcal G}(\Gamma)}(\Delta(\mu, p)(\Gamma), M_1) = 0
\qquad
\text{for}
\quad
(\mu,p)\notin \mathcal S(\lambda, r),
\\
\text{and}
\qquad
\Ext^1_{\widehat{\mathcal G}(\Gamma)}(\Delta(\mu_0, p_0)(\Gamma), M_1) = 0.
\end{gather*}
The fact that $\Ext^1_{\widehat{\mathcal G}(\Gamma)}(\Delta(\mu_1, p_1)(\Gamma), M_1) = 0$ follows from
the fact that $M_1$ is a~summand of $M_1'$ and that $\Ext^1_{\widehat{\mathcal
G}}(\Delta(\mu_1,p_1)(\Gamma),M_1')=0$.

We use Lemma~\ref{elem} again, with $N=M_1$ and $M=\Delta(\mu_2,p_2)(\Gamma)$, and we get
\begin{gather*}
0\to M_1 \to M_2' \to {\Delta(\mu_2,p_2)(\Gamma)}^{d_2'} \to 0,
\end{gather*}
a~summand $M_2\subseteq M_2'$ containing $(M_2')_{\mu_i}[p_i]$, $i=0,1$, such that
\begin{gather*}
M_1 \stackrel{\iota_1}{\hookrightarrow}M_2,
\qquad
\Ext^1_{\widehat{\mathcal G}(\Gamma)}(\Delta(\mu_i, p_i)(\Gamma), M_2) = 0
\qquad
\text{for}
\quad
i=0,1,2,
\\
\text{and}
\qquad
\Ext^1_{\widehat{\mathcal G}(\Gamma)}(\Delta(\mu, p)(\Gamma), M_2) = 0
\qquad
\text{for}
\quad
(\mu,p)\notin \mathcal S(\lambda, r).
\end{gather*}
Repeating this procedure, and using Lemma~\ref{elem} and the properties of $\eta$, we have the following proposition.
The condition on weights is a~consequence that $ \wt \Delta(\lambda_i, p)\subset \conv
W\lambda_k$ for all $i\le k$.

\begin{Proposition}
There exists a~family $\{M_s\}$, $s\in \mathbb Z_{\ge0}$, of indecomposable finite-dimensional modules and injective
morphisms $\iota_s: M_s\to M_{s+1}$ of objects of $\mathcal G_{\bdd}(\Gamma)$ which have the following
properties:
\begin{enumerate}\itemsep=0pt
\item[$1.$]
$M_0=\Delta(\lambda_k,r)(\Gamma)=\Delta(\mu_0,p_0)(\Gamma)$, and for $s\ge 1$,
\begin{gather*}
M_{s}/\iota_{s-1}(M_{s-1})\cong\Delta(\mu_s,p_s)(\Gamma)^{d_s},
\qquad
d_s\in\mathbb Z_+,
\\
\dim M_s[r]_{\lambda_k}=1,
\qquad
\wt M_s\subset\conv W\lambda_k.
\end{gather*}
\item[$2.$]
The spaces $M_s[p]=0$, for all
$s\ge 0$,
$p \gg \max\{r_i \}$.
\item[$3.$]
For all $0\le\ell\le s$ we have $\Ext^1_{\widehat{\mathcal G}}(\Delta(\mu_\ell,p_\ell)(\Gamma), M_s)=0$,
and, for all $(\mu,p)\notin \mathcal S(\lambda, r)$, we have $\Ext^1_{\widehat{\mathcal
G}}(\Delta(\mu,p)(\Gamma), M_s)=0$.
\item[$4.$]
$M_s$ is generated as a~$\mathfrak g[t]$-module by the spaces $\{M_s[p_\ell]_{\mu_\ell}: \ell\le s\}$.
Moreover, if we let
\begin{gather*}
\iota_{r,s}= \iota_{s-1}\cdots\iota_r: \ M_r\to M_{s},
\quad
r< s,
\qquad
\iota_{r,r}=\id,
\end{gather*}
then $M_s[p_\ell]_{\mu_\ell}= \iota_{\ell,s}(M_\ell[p_\ell]_{\mu_\ell})$,
$s\ge\ell$.
\end{enumerate}
\end{Proposition}

\subsection{Def\/ining the tilting modules}\label{section5.3}

Let $T(\lambda_k,r)(\Gamma)=T(\lambda,r)(\Gamma)$ be the direct limit of $\{ M_s, \iota_{r,s} \,|\, r,s \in \mathbb Z_+,
r\le s\}$.
We have an injection $M_s \hookrightarrow T(\lambda,r)(\Gamma)$, and, letting $\widetilde{M_s}$ the image of $M_s$, we
have $\widetilde{M_{s}}\subset \widetilde{M_{s+1}}$, $T(\lambda,r)(\Gamma)= \bigcup {\widetilde{M_s}}$, and
$\frac{\widetilde{M_s}}{\widetilde{M_{s-1}}} \cong \frac{M_s}{M_{s-1}}$.
In particular we see that $T(\lambda,r)(\Gamma)$ has $\Delta(\Gamma)$-f\/iltration.
We identify~$M_s$ with~$\widetilde{M_s}$.

The argument that $T(\lambda,r)(\Gamma)$ is indecomposable is identical to that from~\cite{bc}, which we include for
completeness.
We begin with an easy observation:
\begin{gather}
\label{ms1}
T(\lambda, r)(\Gamma)[p_\ell]_{\mu_\ell} = M_{\ell}[p_\ell]_{\mu_\ell},
\qquad
M_s=\sum\limits_{\ell\le s}\textbf U(\mathfrak g[t]) T(\lambda_k,r)[p_\ell]_{\mu_\ell}.
\end{gather}

To prove that $T(\lambda,r)(\Gamma)$ is indecomposable, suppose that $T(\lambda,r)(\Gamma)= U_1\oplus U_2$.
Since $\dim T(\lambda,r)(\Gamma)[r]_{\lambda}=1$, we may assume without loss of generality that
$T(\lambda,r)(\Gamma)[r]_{\lambda}\subset U_1$ and hence $ M_0\subset U_1$.
Assume that we have proved by induction that $M_{s-1}\subset U_1$.
Since~$M_s$ is gene\-ra\-ted as a~$\mathfrak g[t]$-module by the spaces $\{M_s[p_\ell]_{\mu_\ell}: \ell\le s\}$, it suf\/f\/ices
to prove that \mbox{$M_s[p_s]_{\mu_s}\subset U_1$}.
By~\eqref{ms1}, we have $U_i[p_s]_{\mu_s}\subset M_s$ and hence
\begin{gather*}
M_s= \left (M_{s-1}+ \textbf U(\mathfrak g[t])U_1[p_s]_{\mu_s}\right)\oplus\textbf U(\mathfrak g[t])U_2[p_s]_{\mu_s}.
\end{gather*}
Since $M_s$ is indecomposable by construction, it follows that $U_2[p_s]_{\mu_s}=0$ and $M_s\subset U_1$ which completes
the inductive step.

\begin{Proposition}
Let $(\lambda,r)\in\Gamma$.
\begin{enumerate}\itemsep=0pt
\item[$1.$]
Then there exists an indecomposable module $T(\lambda,r)(\Gamma) \in \Ob {\mathcal
G}_{\bdd}(\Gamma)$ which admits a~filtration by finite-dimensional modules $M_s= \sum\limits_{\ell\le
s}\textbf U(\mathfrak g[t]) T(\lambda,r)(\Gamma)[p_\ell]_{\mu_\ell}$, $s\ge 0$, such that
$M_0\cong\Delta(\lambda,r)(\Gamma)$ and the successive quotients are isomorphic to a~finite-direct sum of
$\Delta(\mu,s)(\Gamma)$, $(\mu,s)\in \mathcal S(\lambda, r)$.
\item[$2.$]
We have $ \wt T(\lambda,r)(\Gamma) \subset\conv W\lambda$,
$\dim T(\lambda,r)(\Gamma)[r]_{\lambda}=1$.
\item[$3.$]
For all $(\mu,s)\in \Gamma$, we have $\Ext^1_{\widehat{\mathcal G}(\Gamma)}(\Delta(\mu,s)(\Gamma)),
T(\lambda,r)(\Gamma))=0$.
\end{enumerate}
\end{Proposition}

\begin{proof} Part (1) and (2) are proved in the proceeding discussion.
The proof for part (3) is identical to that found in~\cite{bc}.
\end{proof}

\subsection{Initial properties of tilting modules} \label{section5.4}

The next result is an analog of Fitting's lemma for the
inf\/inite-dimensional modules $T(\lambda,r)(\Gamma)$.
\begin{Lemma}
Let $\psi: T(\lambda,r)(\Gamma)\to T(\lambda,r)(\Gamma)$ be any morphism of objects of $\widehat{\mathcal G}$.
Then $\psi(M_s)\subset M_s$ for all $s\ge 0$ and $\psi$ is either an isomorphism or locally nilpotent, i.e., given $m\in
M$, there exists $\ell\ge 0$ $($depending on $m)$ such that $\psi^\ell(m)=0$.
\end{Lemma}

\begin{proof} Since $\psi$ preserves both weight spaces and graded components it follows that \mbox{$\psi(M_s)\!\subset\! M_s$} for all
$s\ge 0$.
Moreover, since $M_s$ is indecomposable and f\/inite--dimensional it follows from Fitting's lemma that the restriction of
$\psi$ to $M_s$, $s\ge 0$ is either nilpotent or an isomorphism.
If all the restrictions are isomorphisms then since $T(\lambda,r)(\Gamma)$ is the union of $M_s$, $s\ge 0$, it follows
that $\psi$ is an isomorphism.
On the other hand, if the restriction of $\psi$ to some $M_s$ is nilpotent, then the restriction of $\psi$ to all
$M_\ell$, $\ell\ge 0$ is nilpotent which proves that $\psi$ is locally nilpotent.
\end{proof}

In the rest of the section we shall complete the proof of the main theorem by showing that any indecomposable tilting
module is isomorphic to some $T(\lambda,r)(\Gamma)$ and that any tilting module in $\mathcal
G_{\bdd}(\Gamma)$ is isomorphic to a~direct sum of indecomposable tilting modules.
Let $T\in\mathcal G_{\bdd}(\Gamma)$ be a~f\/ixed tilting module.
Then we have
\begin{gather}
\label{tiltext}
\Ext^1_{\widehat{\mathcal G}}(T,\nabla(\lambda,r)(\Gamma))=0=\Ext^1_{\widehat{\mathcal
G}}(\Delta(\lambda,r)(\Gamma), T),
\qquad
(\lambda,r)\in \Gamma,
\end{gather}
where the f\/irst equality is due to Corollary~\ref{extMnabla}.

\begin{Lemma}
Suppose that $T_1$ is any summand of $T$.
Then, $T_1$ admits a~$\nabla(\Gamma)$-filtration and $\Ext^1_{\widehat{\mathcal
G}}(T_1,\nabla(\lambda,r)(\Gamma))=0$, for all $(\lambda,r)\in\Gamma$.
\end{Lemma}
\begin{proof} Since $\Ext^1$ commutes with f\/inite direct sums, for all $ (\lambda,r)\in \Gamma$ we have
\begin{gather*}
\Ext^1_{\widehat{\mathcal G}}(T_1,\nabla(\lambda,r) (\Gamma))=0,
\qquad
\text{and}
\qquad
\Ext^1_{\widehat{\mathcal G}}(\Delta(\lambda,r)( \Gamma), T_1)=0.
\end{gather*}
By Proposition~\ref{extnablaconn}, the second equality implies that $T_1$ has a~$\nabla (\Gamma)$-f\/iltration and the
proof of the lemma is complete.
\end{proof}

\subsection{Completing the proof of the main theorem} \label{section5.5}

The preceding lemma illustrates one of the dif\/f\/iculties we face in
our situation.
Namely, we cannot directly conclude that $T_1$ has a~$\Delta( \Gamma)$-f\/iltration from the vanishing
$\Ext$-condition by using the dual of Proposition~\ref{extnablaconn}.
However, we have the following, whose proof is given in~\cite{bc}.

\begin{Proposition}
\label{summands}
Suppose that $N\in\mathcal G_{\bdd} (\Gamma)$ has a~$\nabla (\Gamma)$-filtration and satisfies
\begin{gather*}
\Ext^1_{\widehat{\mathcal G}}(N,\nabla(\lambda,r) (\Gamma)) =0,
\qquad
\text{for all}
\quad
(\lambda,r)\in \Gamma.
\end{gather*}
Let $(\mu,s)$ be such that $N\rightarrow \nabla(\mu,s)(\Gamma)\rightarrow 0$.
Then $T(\mu,s)(\Gamma)$ is a~summand of $N$.
\end{Proposition}

The following is immediate.
Note that this also shows that our construction of the indecomposable tilting modules is independent of the choice of
enumeration of $P^+$, the set $\mathcal S(\lambda, r)$ and~$\eta$.

\begin{Corollary}
Any indecomposable tilting module is isomorphic to $T(\lambda{,}r)(\Gamma)$ for some \mbox{$(\lambda{,}r)\!\in\! \Gamma.\!$}
Further if $T\in\Ob{\mathcal G}_{\bdd}(\Gamma)$ is tilting there exists $(\lambda,r)\in
\Gamma$ such that $T(\lambda,r)(\Gamma)$ is isomorphic to a~direct summand of $T$.
\end{Corollary}
\begin{proof} Since $T$ and $T(\lambda,r)(\Gamma)$ are tilting they satisfy~\eqref{tiltext} and the corollary follows.
\end{proof}

We can now prove the following theorem.

\begin{Theorem}
Let $T\in\Ob{\mathcal G}_{\bdd}(\Gamma)$.
The following are equivalent.
\begin{enumerate}\itemsep=0pt
\item[$1.$]
$T$ is tilting.
\item[$2.$]
$\Ext^1_{\widehat{\mathcal G}}(\Delta(\lambda,r)(\Gamma), T) =0 =\Ext^1_{\widehat{\mathcal
G}}(T,\nabla(\lambda,r)(\Gamma))$,
$(\lambda,r)\in \Gamma$.
\item[$3.$]
$T$ is isomorphic to a~direct sum of objects $T(\mu,s)(\Gamma)$, $(\mu,s)\in \Gamma $.
\end{enumerate}
\end{Theorem}

\begin{proof} The implication (1)$\implies$(2) is given by Corollary~\ref{extMnabla} and Proposition~\ref{extnablaconn}, while
the fact that (3) implies (1) is clear.
We complete the proof by showing that (2) implies (3).

By Proposition~\ref{26}, that $\Ext^1_{\widehat{\mathcal G}}(\Delta(\lambda,r)(\Gamma), T) =0 $ implies
that $T$ admits a~$\nabla(\Gamma)$-f\/iltration.
By Lemma~\ref{reorder}, we can assume that the f\/iltration is f\/inite, and if $\lambda_k=\lambda$ is maximal such that
$T_\lambda\ne0$, then $T\rightarrow \nabla(\lambda, r)(\Gamma)\rightarrow 0$ for $r$.
Indeed, if we choose integers $r_1\ge r_2\ge \cdots $ such that
\begin{gather*}
[T:\nabla(\lambda,s)(\Gamma)]\ne 0
\qquad
\text{if\/f}
\quad
s=r_i
\quad
\text{for some}
\quad
i,
\end{gather*}
then Lemma~\ref{reorder} says that we may assume that $T\rightarrow \nabla(\lambda, r_1)(\Gamma)\rightarrow 0$.
By Proposition~\ref{summands} we have $T\cong T(\lambda, r_1)(\Gamma)\oplus T_1$ and we see that $T_1$ has
a~$\nabla(\Gamma)$-f\/iltration.
The same argument implies that $T_1$ maps onto $\nabla(\lambda, r_2)(\Gamma)$, and hence $T_1 \cong T(\lambda,
r_2)(\Gamma) \oplus T_2$.

Continuing, we f\/ind that for $j\ge 1$, there exists a~summand $T_j$ of $T$ with
\begin{gather*}
T=T_j \bigoplus_{s=1}^j T(\lambda, r_s)(\Gamma).
\end{gather*}
Let $\pi_j:T\rightarrow \oplus_{s=1}^j T(\lambda, r_s)(\Gamma)$ be the canonical projections.
Because $T$ has f\/inite-dimensional graded components, and the $r_i$ are decreasing, it follows that for $m\in T$ there
exists an integer~$k(m)$ such that $\pi_j(m)=\pi_{k(m)}(m)$ for all $k(m)\le j$.
Hence, we have a~surjection
\begin{gather*}
\pi: \ T\rightarrow \bigoplus_{j\ge 1} T(\lambda, r_j)(\Gamma) \rightarrow 0
\qquad
\text{and}
\qquad
\ker \pi = \bigcap T_j,
\end{gather*}
where $\pi(m):=\pi_{k(m)}(m)$.
In particular, we have $T=(\bigoplus T(\lambda, r_i)(\Gamma)) \oplus \ker \pi$, where $(\ker \pi)_\lambda=0$, $\ker \pi$
admits a~$\nabla(\Gamma)$-f\/iltration and $\Ext^1_{\widehat{\mathcal G}}(\ker
\pi,\nabla(\mu,r)(\Gamma))=0$, for all $(\mu,r)\in \Gamma$.
It follows that we may apply to $\ker \pi$ the same arguments we used on~$T$.
The result follows by induction on~$k$.
\end{proof}

\subsection{Proof of Proposition~\ref{subset}}\label{section5.6}

We construct here the set $\mathcal S(\lambda, r)$ and the enumeration $\eta$.

\textit{The set $\mathcal S(\lambda, r)$.} Recall that $\Gamma=P^+\times J$, and that $a=\inf  J$ and
$b=\sup  J$.
Using the enumeration of $P^+$, let $\lambda=\lambda_k$ and def\/ine integers $r_k \le r_{k-1} \le \dots \le r_0$
recursively by setting $r_k = r$ and
\begin{gather*}
r_s = \max \{ r \,|\, \Delta(\lambda_{s+1},r_{s+1})(\Gamma)[r] \ne 0\}.
\end{gather*}
Note that because $\Delta(\lambda_i, r_i)(\Gamma)[p]=0$ for any $p> r_{i-1}$, and $r_i\le r_j$ for all $j<i$, we have
$\Delta(\lambda_i, r_i)(\Gamma)[p]=0$ if $p>r_j$ for any $j<i$.
Then, it follows that
\begin{gather}
\label{eq4.6}
\Delta(\lambda_i, s)(\Gamma)[p]=0
\qquad
\text{for any}
\quad
s \le r_i < p.
\end{gather}

We set $\mathcal S(\lambda, r)=\{ (\lambda_i,s)| i\le k, s\le r_i \}$, and note that it satisf\/ies conditions (1) and
(2) of Proposition~\ref{subset} by construction.

We now verify condition (3).
Let $(\mu,s)\in \mathcal S(\lambda, r)$ and $(\mu',s')\notin \mathcal S(\lambda, r)$.
There are two possibilities for $(\mu',s')$: either $\mu'=\lambda_i$ for $i> k$, or $\mu'=\lambda_i$ for some $0\le i
\le k$ and $s'>r_i$.
The f\/irst case is covered by Proposition~\ref{homresults}.3, which tells us that if $\lambda \not\le \mu$ then
$\Ext^1_{\widehat{\mathcal G}}(\Delta(\lambda, r)(\Gamma), \Delta(\mu,s)(\Gamma))=0$ for all $s,r\in
\mathbb Z$.
For the second case, again using Proposition~\ref{homresults}.3, it is enough to prove that
$\Ext^1_{\widehat{\mathcal G}}(\Delta(\lambda_i, s')(\Gamma), \Delta(\lambda_j, s)(\Gamma))=0$ for
$\lambda_i \le \lambda_j$, $s'>r_i$, and $s\le r_j$.
By the total order (cf.\
Section~\ref{section4.4}), it follows that $i\le j$, and, hence, $s\le r_j \le r_i < s'$.
By Proposition~\ref{homresults}.4, we can in fact assume that $i < j$.
Consider a~short exact sequence $0 \to \Delta(\lambda_j, s)(\Gamma) \to M \to \Delta(\lambda_i, s')(\Gamma) \to 0$.
From~\eqref{eq4.6} we have $\Delta(\lambda_j, s)(\Gamma)[s']=0$, and it follows that the sequence splits, as we desired.

\textit{The enumeration $\eta$.} It remains to def\/ine the enumeration $\eta$.
The case where $J=\mathbb Z$ is done in~\cite{bc}, and the case where $J$ is a~f\/inite or of the form $[a, \infty)$ (in
which case $\mathcal S(\lambda, r)$ is in fact f\/inite), we use the enumeration def\/ined by the following rules
\begin{enumerate}\itemsep=0pt
\item[1)]
$\eta(\lambda_i, s) < \eta(\lambda_j, s')$ if $i>j$,
\item[2)]
$\eta(\lambda_i, s) < \eta(\lambda_i, s-1)$.
\end{enumerate}
Suppose that $i<j$ and let $(\mu_i,p_i)=\eta^{-1}(i)$ and $(\mu_j,p_j)=\eta^{-1}(j)$.
If $\mu_i=\mu_j$, then rule~$(2)$ implies that $p_j < p_i$, and Proposition~\ref{homresults} says that
$\Ext^1_{\widehat{\mathcal G}}(\Delta(\mu_i,p_i)(\Gamma),\Delta\mu_j,p_j)(\Gamma))=0$.
Otherwise, we have $\mu_j \not\le \mu_i$, and the result again follows by Proposition~\ref{homresults}.

We are left with the case where $J=(-\infty, b]$.
In this case $\eta$ will in fact be a~bijection.
Note that it is enough to def\/ine a~bijective, set theoretic inverse $\eta^{-1}$.

We recursively def\/ine another set of integers $\{ r_i'\}$ by setting $r_k'=r_k$ and letting $r_i'=\max\{ r|
\Delta(\lambda_{i+1}, r_{i+1}')[r]\ne 0\}$.
It is easy to see that $r_i \le r_i'$.
If $r_i < r_i'$, then we must have $r_i'> b$, which implies that $r_i=b$.
This implies that $r_j=b$ for all $j<i$.
We note that $\Delta(\lambda_j,b)(\Gamma)=V(\lambda_j, b)$.
Set $a_s:= r_s'-r_{s+1}'$.

\begin{Lemma}
\label{thelaststop}
We have $\Ext^1_{\widehat{\mathcal G}}(\Delta(\lambda_i,c)(\Gamma),\Delta(\lambda_{s},d)(\Gamma))=0$
if $c-d\ge a_{s-1}+1$ and $i<s$.
\end{Lemma}
\begin{proof} We f\/irst prove that under these conditions $\Ext^1_{\widehat{\mathcal
G}}(\Delta(\lambda_i,c),\Delta(\lambda_{s},d))=0$.
Note that we can shift by $-d+r_s'$, and so we examine $\Ext^1_{\widehat{\mathcal
G}}(\Delta(\lambda_i,c-d+r_s'),\Delta(\lambda_{s},r_s'))$.
According to our hypothesis, we have $c-d+r_s' \ge r_{s-1}'+1$.
It follows by the def\/inition of $r_{s-1}'$ that $\Delta(\lambda_s,r_s')[p]=0$ if $p\ge c-d+r_s'$.
If we examine a~sequence $0 \to \Delta(\lambda_s, r_s') \to M \to \Delta(\lambda_i, c-d+r_s') \to 0$ and let $m\in
M_{\lambda_i}[c-d+r_s']$ be a~pre-image of $w_{\lambda_i}$, it is clear that $m$ satisf\/ies the def\/ining relations of
$w_{\lambda_i}$.
Therefore, the sequence splits.
Finally, observe that $\Delta(\lambda_s, d)(\Gamma)[c]=0$, again by the def\/inition of the $r_i'$ and noting that
$\Delta(\lambda_s, d)(\Gamma) $ is a~quotient of $\Delta(\lambda_s, d)$.
The same argument as above also shows that $\Ext^1_{\widehat{\mathcal
G}}(\Delta(\lambda_i,c)(\Gamma),\Delta(\lambda_{s},d)(\Gamma))=0$.
\end{proof}

We def\/ine $\eta^{-1}: \mathbb Z_{\ge 0} \to \mathcal S(\lambda, r)$ in the following way.

Set $\eta^{-1}(0) = (\lambda_k,r_k)$.
If $\eta^{-1}$ is def\/ined on $\{0, \ldots, m-1\}$ and $\eta^{-1}(m-1)=(\lambda_i, p_i)$, we def\/ine $\eta^{-1}(m)$ as
follows.
Suppose that $i>0$, and $(\lambda_{i-1}, p_i + a_{i-1})\in \mathcal S(\lambda, r)$, then we set
$\eta^{-1}(m)=(\lambda_{i-1}, p_i + a_{i-1})$.
Otherwise, we let $\eta^{-1}(m)=(\lambda_k, p_m-1)$, where $p_m$ is the minimal integer such that $(\lambda_k, p)$ has
a~pre-image under $\eta^{-1}$.

To prove that $\Ext^1_{\widehat{\mathcal G}}(\Delta(\mu_i,p_i)(\Gamma),\Delta(\mu_j,p_j)(\Gamma))=0$ if
$i< j$, we assume that $\mu_i \le \mu_j$.
If $\mu_i = \mu_j$ then this follows from Proposition~\ref{homresults}.
So we assume that $\mu_i < \mu_j$, and lets say that $\mu_j=\lambda_\ell$.
In this case we must have $p_i-p_j > a_{\ell-1}$ and so the result follows by Lemma~\ref{thelaststop}.

\section{Some dif\/ferent considerations on truncated categories}\label{section6}

Throughout this section we discuss some ``trivial'' tilting theories for the category $\mathcal G(\Gamma)$ by
considering dif\/ferent type of orders on the set $\Lambda$.
These categories, equipped with the orders described below, have already appeared in the literature
(see~\cite{BCFM,cg:hwcat,cg:koszul,cg:minp,ckr:faces} and references therein).

\subsection{Truncated categories with the covering relation}\label{section6.1}

Consider a strict partial order on $\Lambda$ in the following way. Given $(\lambda,r), (\mu,s)\in\Lambda$, we say that
\begin{gather*}
(\mu,s) \text{ covers } (\lambda,r) \text{ if and only if } s=r+1 \text{ and } \mu-\lambda\in R\cup\{0\}.
\end{gather*}
Notice that for any $(\mu,s)\in\Lambda$ the set of $(\lambda,r)\in\Lambda$ such that $(\mu,s)$ covers $(\lambda,r)$ is f\/inite. Let $\preceq$ be the unique partial order on~$\Lambda$ generated by this covering relation.

One of the main inspirations to consider this relation comes from the following proposition:

\begin{Proposition}[\protect{\cite[Proposition 2.5]{cg:hwcat}}] For $(\lambda,r),(\mu,s)\in\Lambda$, we have
\begin{gather*}
\Ext^1_{{\mathcal G}}(V(\lambda,r),V(\mu,s))=
\begin{cases}
0, &\text{if}\quad s\ne r+1,
\\
\Hom_{\mathfrak g}(V(\lambda), \mathfrak g\otimes V(\mu)), &\text{if}\quad s= r+1.
\end{cases}
\end{gather*}
In other words, $\Ext^1_{{\mathcal G}}(V(\lambda,r),V(\mu,s))=0$ except when $(\mu,s)$ covers $(\lambda,r)$.
\end{Proposition}

Given $\Gamma\subset\Lambda$, set
\begin{gather*}
{V_\Gamma}^+= \{v\in V[s]_\mu: \mathfrak n^+ v=0,
\
(\mu,s)\in \Gamma\},
\qquad
V_\Gamma=\textbf U(\mathfrak g)V_\Gamma{}^+,
\qquad
V^\Gamma=V/V_{\Lambda\setminus\Gamma}.
\end{gather*}

\begin{Proposition}[\protect{\cite[Propositions 2.1, 2.4, and 2.7]{cg:hwcat}}]
\label{propofPandI}
Let $\Gamma$ be finite and convex and assume that $(\lambda,r), (\mu,s)\in\Gamma$.
\begin{enumerate}\itemsep=0pt
\item[$1.$]
$[P(\lambda,r)^\Gamma:V(\mu,s)]=[P(\lambda,r):V(\mu,s)]=[I(\mu,s):V(\lambda,r)]=[I(\mu,s)_\Gamma:V(\lambda,r)]$.

\item[$2.$]
$\Hom_{{\mathcal G}}(P(\mu,s),P(\lambda,r)) \cong
\Hom_{\mathcal{G}[\Gamma]}(P(\mu,s)^{\Gamma},P(\lambda,r)^{\Gamma})$.

\item[$3.$]
Let $K(\lambda,r )$ be the kernel of the canonical projection $P(\lambda,r)\twoheadrightarrow V(\lambda,r)$, and let
$(\mu,s)\in\Lambda$.
Then $[K(\lambda,r):V(\mu,s)]\not=0$ only if $(\lambda,r)\prec(\mu,s)$.

\item[$4.$]
Let $(\mu,s)\in\Lambda$.
Then~$[I(\lambda,r)/V(\lambda,r):V(\mu,s)]\not=0$ only if~$(\mu,s)\prec(\lambda,r)$.
\end{enumerate}
\end{Proposition}

Following Section~\ref{section2.6} but using the partial order $\preceq$ def\/ined above, for each $(\lambda,r)\in \Gamma$
we denote by $\Delta(\lambda,r)(\Gamma)$ the maximal quotient of $P(\lambda,r)$ such that
\begin{gather*}
[\Delta(\lambda,r)(\Gamma):V(\mu,s)] \ne 0\implies (\mu,s)\preceq (\lambda,r).
\end{gather*}
Similarly, we denoted by $\nabla(\lambda,r)(\Gamma)$ the maximal submodule of $I(\lambda,r)$ satisfying
\begin{gather*}
[\nabla(\lambda,r)(\Gamma):V(\mu,s)] \ne 0\implies (\mu,s)\preceq (\lambda,r).
\end{gather*}
The modules $\Delta(\lambda,r)(\Gamma)$ and $\nabla(\lambda,r)(\Gamma)$ are called, respectively, the standard and
co-standard modules associated to $(\lambda,r)$.
Further, any module in $\widehat{\mathcal G}$ with a~$\Delta(\Gamma)$-f\/iltration and a~$\nabla(\Gamma)$-f\/iltration is
called tilting.

\begin{Proposition}\label{p:tilt}
Let $\Gamma \subseteq \Lambda$ finite and convex.
\begin{enumerate}\itemsep=0pt
\item[$1.$]
For all $(\lambda,r)\in \Gamma$, there exists indecomposable tilting module $T(\lambda,r)(\Gamma)$ and
$T(\lambda,r)(\Gamma)=I(\lambda,r)$.
\item[$2.$]
For all indecomposable tilting module $T$, we have $T\cong T(\lambda,r)(\Gamma)$ for some $(\lambda,r) \in \Gamma$.
\item[$3.$]
Every tilting module is isomorphic to a~direct sum of indecomposable tilting modules.
\end{enumerate}
\end{Proposition}

Before proving this proposition, we have some remarks to make:

\begin{Remark}\label{stdcostd}\quad
\begin{enumerate}\itemsep=0pt
\item
It follows from Proposition~\ref{projectives}.3 and Proposition~\ref{propofPandI},
parts~(1) and~(4), that the costandard modules in $\widehat{\mathcal G}$ associated to
$(\lambda,r)$ is $I(\lambda, r)$ and similarly it follows from Proposition~\ref{propofPandI},
parts~(1) and~(3), that the standard module in $\widehat{\mathcal G}$ associated to
$(\lambda,r)$ is the simple module $V(\lambda,r)$.

\item
For any $M\in \mathcal G$ let $k(M)$ the such that $M\in \mathcal G_{\leq k(M)}$.
Thus $M$ admits a~f\/iltration $\{ M_i\} $ where $M_i = \bigoplus_{j=0}^i M[k(M)-j]$ which can be ref\/ined into
a~Jordan-Holder series since each quotient $M_{i+1}/M_i$ is a~f\/inite-dimensional $\mathfrak g$-module.
\end{enumerate}
\end{Remark}

\begin{proof}
[Proof of Proposition~\ref{p:tilt}] Part (1) follows from the Remark~\ref{stdcostd} since we have $T(\lambda,
r)(\Gamma):= I(\lambda,r)$.
Part (2) and (3) are direct consequences of the injectivity of $I(\lambda,r)$.
\end{proof}

\subsection{Truncated categories related to restricted Kirillov--Reshetikhin modules}\label{section6.2}

One of the goals of~\cite{BCFM,cg:minp} was to study the modules $P(\lambda, r)^\Gamma$ (and their multigraded version)
under certain very specif\/ic conditions on $\Gamma$.
In these papers it was shown that the modules $P(\lambda,r)^\Gamma$ are giving in terms of generators and relations
which allows us to regard these modules as specializations of the famous Kirillov--Reshetikhin modules (in the sense
of~\cite{CM:kr,CM:krg}).
These papers develop a~general theory over a~$\mathbb Z_+$-graded Lie algebra $\mathfrak a~= \bigoplus_{i\in \mathbb Z}
\mathfrak g[i]$ where $\mathfrak g_0$ is a~f\/inite-dimensional complex simple Lie algebra and its non-zero graded
components $\mathfrak g[i]$, $i>0$, are f\/inite-dimensional $\mathfrak g_0$- modules.
By focusing in these algebras with $\mathfrak g[i] =0$ for $i>1$, we have $\mathfrak a~\cong \mathfrak g\ltimes V$,
where $V$ is a~$\mathfrak g$-module, and in this context a~very particular tilting theory can be described as follows.

Assume that $\wt(V)\ne\{0\}$ and f\/ix a~subset $\Psi\subseteq\wt(V)$ satisfying
\begin{gather*}
\sum\limits_{\nu\in\Psi} m_\nu\nu= \sum\limits_{\mu\in\wt(V)} n_\mu\mu(m_\nu,n_\mu\in\mathbb Z_+)
\
\Longrightarrow
\
\sum\limits_{\nu\in\Psi} m_\nu\le \sum\limits_{\mu\in\wt(V)} n_\mu
\end{gather*}
and
\begin{gather*}
\sum\limits_{\nu\in\Psi} m_\nu= \sum\limits_{\mu\in\wt(V)}n_\mu
\qquad
\text{only if}
\quad
n_\mu=0
\quad
\text{for all}
\quad
\mu\notin\Psi.
\end{gather*}

\begin{Remark}
Such subsets are precisely those contained in a proper face of the convex polytope determined by $\wt (V)$ conform~\cite{kharid:polytopes}.
\end{Remark}

Consider the ref\/lexive and transitive binary relation on $P$ given by
\begin{gather*}
\mu\le_\Psi\lambda
\qquad
\text{if}
\quad
\lambda-\mu\in\mathbb Z_+\Psi,
\end{gather*}
where $\mathbb Z_+\Psi$ is the $\mathbb Z_+$-span of $\Psi$.
Set also
\begin{gather*}
d_\Psi(\mu,\lambda)=\min\left\{\sum\limits_{\nu\in\Psi}m_\nu: \lambda-\mu=\sum\limits_{\nu\in\Psi}m_\nu\nu,\,
m_\nu\in\mathbb Z_+ \;
\forall\,
\nu\in\Psi \right\}.
\end{gather*}
By \cite[Proposition~5.2]{ckr:faces}, $\le_\Psi$ is in fact a partial order on~$P$.
Moreover, it induces a ref\/inement $\preccurlyeq_\Psi$ of the partial order~$\preceq$ on~$\Lambda$ by setting
\begin{gather*}
(\lambda, r)\preccurlyeq_\Psi (\mu, s)
\qquad
\text{if}
\quad
\lambda\le_\Psi\mu,
\quad
s- r\in\mathbb Z_+,
\quad
\text{and}
\quad
d_\Psi(\lambda,\mu)= s- r.
\end{gather*}

Finally, if $\Gamma\subseteq\Lambda$ is f\/inite and convex with respect to $\preccurlyeq_\Psi$ and there exists
$(\lambda, r)\in\Lambda$ such that $(\lambda, r)\preccurlyeq_\Psi(\mu, s)$ for all $(\mu, s)\in\Gamma$, it was
shown~\cite[Lemma 5.5]{ckr:faces} that
\begin{gather}
\label{hom}
\Hom_{\mathcal G[\Gamma]}(P(\mu, s)^\Gamma,P(\nu, t)^\Gamma)\ne 0
\qquad
\text{only if}
\quad
(\nu, t)\preccurlyeq_\Psi (\mu, s).
\end{gather}
In particular, it follows from Proposition~\ref{projectives}.3, Proposition~\ref{propofPandI},
parts~(2) and~(3), and~\eqref{hom} that
\begin{gather*}
[P(\lambda,r):V(\mu,s)]\ne 0 \implies (\lambda,r)\preccurlyeq_\Psi (\mu,s)
\end{gather*}
and
\begin{gather*}
[I(\lambda,r):V(\mu,s)]\ne 0 \implies (\mu,s)\preccurlyeq_\Psi (\lambda,r).
\end{gather*}
We conclude that the standard modules in $\mathcal G(\Gamma)$ are the simple modules and the costandard modules are the
injective hulls of the simple modules and, hence, $T(\lambda,r)(\Gamma)=I(\lambda,r)_\Gamma$.

\subsection*{Index of notation}

\hspace{\parindent}Subsection~\ref{section1.2} --- $\mathcal F(\mathfrak g)$, $V(\lambda)$.

Subsection~\ref{section2.1} --- $\widehat{\mathcal G}$, $V[r]$, $V(\lambda,r)$, $v_{\lambda,r}$, $M^*$, $\Lambda$, $\le$.

Subsection~\ref{section2.2} --- $\mathcal G_{\le s}$, $\mathcal G$, $\mathcal G_{\bdd}$, $V_{>s}$, $V_{\le
s}$, $[V:V(\lambda,s)]$, $\Lambda(V)$.

Subsection~\ref{section2.3} --- $P(\lambda,r)$, $I(\lambda,r)$, $p_{\lambda,r}$.

Subsection~\ref{section2.4} --- $\Delta(\lambda,r)$, $W(\lambda,r)$, $\nabla(\lambda,r)$.

Subsection~\ref{section2.5} --- $\Gamma$, $\widehat{\mathcal G}(\Gamma)$, $\mathcal G(\Gamma)$, $\mathcal G_{\bdd}(\Gamma)$.

Subsection~\ref{section2.6} --- $J$, $a$, $b$, $M^\Gamma$, $\Delta(\lambda,r)(\Gamma)$, $W(\lambda,r)(\Gamma)$, $\nabla(\lambda,r)(\Gamma)$.

Section~\ref{section3} --- $\Delta(\Gamma)$ (respect.
$\nabla(\Gamma)$) -f\/iltration$, $$T(\lambda,r)$.

Subsection~\ref{section4.3} --- $I(M)$.

Subsection~\ref{section4.4} --- o-canonical f\/iltration.

Subsection~\ref{section4.5} --- $N_s$, $I(\lambda,r)_p$.

Subsection~\ref{section4.6} --- $T(\lambda,r)(\Gamma)$.

\subsection*{Acknowledgements}

The authors are grateful for the stimulating discussions with Professor Adriano Moura as well as the hospitality of the
Institute of Mathematics of the State University of Campinas where most of this work was completed.
We thank the anonymous referees for their signif\/icant contributions towards improving the exposition of this paper.
This work was partially supported by FAPESP grants 2012/06923-0 (M.~Bennett) and 2011/22322-4 (A.~Bianchi).

\pdfbookmark[1]{References}{ref}
\LastPageEnding

\end{document}